\documentclass{amsart}[12pt]
\def\doctype{}

\usepackage{latexsym,amssymb}
\usepackage{color}
\usepackage{dsfont}
\usepackage{fancyhdr}
\usepackage{kbordermatrix}
\usepackage{tikz}
\usepackage{hyperref}
\hypersetup{
colorlinks=true,
citecolor=blue,
linkcolor=blue,
  urlcolor=blue
}

\newcommand\one{\mathds{1}}

\newcommand\vb{\mathbf{b}}
\newcommand\vd{\mathbf{d}}
\newcommand\vf{\mathbf{f}}
\newcommand\vx{\mathbf{x}}
\newcommand\vy{\mathbf{y}}

\newcommand\Z{\mathbb{Z}}

\newcommand\R{\mathbb{R}}
\newcommand\tr{\mathrm{tr}}
\newcommand\Aut{\mathrm{Aut}}
\newcommand\bal{\mathrm{bal}}
\newcommand\Sym{\mathcal{S}}
\newcommand\lam{\lambda}
\newcommand\lcm{\mathrm{lcm}}
\def\dotcup{\,\dot\cup\,}

\newcommand{\comment}[1]{}

\numberwithin{equation}{section}

\fancypagestyle{titlepage}{
\fancyhead[R]{\doctype}
\fancyhead[CL]{}
\cfoot{\vspace{5pt} \thepage}
}

\let\oldsection\section
\newcommand\boldsection[1]{\oldsection{\bf #1}}
\newcommand\starsection[1]{\oldsection*{\bf #1}}
\makeatletter
\renewcommand\section{\@ifstar\starsection\boldsection}
\makeatother

\newtheoremstyle{theorem}
  {12pt}		  % space above
  {0pt}  % space below
  {\sl}  % bofy font
  {\parindent}     % ident - empty=no indent,  \parindent= paragraph indent
  {\bf}  % thm head font
  {. }    % punctuation after thm head
  { }    % space after thm head: `` ``=normal \newline=linebreak
  {}     % thm head specification
\theoremstyle{theorem}
\newtheorem{thm}{Theorem}[section]  % 1st argument is your name for it
\newtheorem{lemma}[thm]{Lemma}     % 2nd argument is what is printed
\newtheorem{cor}[thm]{Corollary}

\newtheorem{prop}[thm]{Proposition}

\newtheoremstyle{definition}
  {12pt}		  % space above
  {0pt}  % space below
  {}  % bofy font
  {\parindent}     % ident - empty=no indent,  \parindent= paragraph indent
  {\bf}  % thm head font
  {. }    % punctuation after thm head
  { }    % space after thm head: `` ``=normal \newline=linebreak
  {}     % thm head specification
\theoremstyle{definition}
\newtheorem{defn}[thm]{Definition}
\newtheorem{ex}[thm]{Example}

\newcommand\rk{{\sc Remark.} }
\newcommand\rks{{\sc Remarks.} }

\renewcommand{\proofname}{Proof}

\makeatletter
\renewenvironment{proof}[1][\proofname]{\par
  \pushQED{\qed}%
  \normalfont \partopsep=\z@skip \topsep=\z@skip
  \trivlist
  \item[\hskip\labelsep
        \scshape
    #1\@addpunct{.}]\ignorespaces
}{%
  \popQED\endtrivlist\@endpefalse
}
\makeatother

\makeatletter
\renewcommand*\@maketitle{%
  \normalfont\normalsize
  \@adminfootnotes
  \@mkboth{\@nx\shortauthors}{\@nx\shorttitle}%
  \global\topskip42\p@\relax % 5.5pc   "   "   "     "     "
  \@settitle
  \ifx\@empty\authors \else {\vskip 1em
\vtop{\centering\shortauthors\@@par}} \fi
  \ifx\@empty\@date \else {\vskip 1em \vtop{\centering\@date\@@par}}\fi % MY CHANGE
  \ifx\@empty\@dedicatory
  \else
    \baselineskip18\p@
    \vtop{\centering{\footnotesize\itshape\@dedicatory\@@par}%
      \global\dimen@i\prevdepth}\prevdepth\dimen@i
  \fi
  \@setabstract
  \normalsize
  \if@titlepage
    \newpage
  \else
    \dimen@34\p@ \advance\dimen@-\baselineskip
    \vskip\dimen@\relax
  \fi
} % end \@maketitle
\renewcommand*\@adminfootnotes{%
  \let\@makefnmark\relax  \let\@thefnmark\relax
  \ifx\@empty\@subjclass\else \@footnotetext{\@setsubjclass}\fi
  \ifx\@empty\@keywords\else \@footnotetext{\@setkeywords}\fi
  \ifx\@empty\thankses\else \@footnotetext{%
    \def\par{\let\par\@par}\@setthanks}%
  \fi
\thispagestyle{titlepage}
}
\makeatother

\begin{document}

\title[Balancing permuted copies]{\large Balancing permuted copies of multigraphs\\ and integer matrices}

\author{Coen del Valle}
\address{
School of Mathematics and Statistics,
University of St Andrews,
St Andrews, UK
}
\email{cdv1@st-andrews.ac.uk}

\author{Peter J.~Dukes}
\address{
Department of Mathematics and Statistics,
University of Victoria, Victoria, Canada
}
\email{dukes@uvic.ca}

\date{\today}

\begin{abstract}
Given a square matrix $A$ over the integers, we consider the $\Z$-module $M_A$ generated by the set of all matrices that are permutation-similar to $A$.  Motivated by analogous problems on signed graph decompositions and block designs, we are interested in the completely symmetric matrices $a I + b J$ belonging to $M_A$.  We give a relatively fast method to compute a generator for such matrices, avoiding the need for a very large canonical form over $\Z$.  Several special cases are considered.  In particular, the problem for symmetric matrices answers a question of Cameron and Cioab\v{a} on determining the eventual period for integers $\lam$ such that the $\lam$-fold complete graph $\lam K_n$ has an edge-decomposition into a given (multi)graph.
\end{abstract}

\maketitle
\hrule

\bigskip

\section{Introduction}
\label{sec-intro}

Given an $n \times n$ integer matrix $A$, consider the $\Z$-module $M_A$ generated by all $n!$ conjugates $P^\top A P$ of $A$ by a permutation matrix $P$.  It is easy to see that the sum of all such copies of $A$ equals $a I + b J$, where $I$ and $J$ denote the identity and all-ones matrices, respectively, and the coefficients $a,b$ satisfy 
\begin{align*}
a+b & = (n-1)! \; \mathrm{tr}(A),\\
a+nb &= (n-1)! \;  \mathrm{sum}(A).
\end{align*}
A matrix of the form $a I + b J$ is sometimes called `completely symmetric'.   We are interested here in the problem of classifying the completely symmetric matrices belonging to $M_A$ for a given $A \in \Z^{n \times n}$.

In the special case when $A$ is a (nonzero) $\{0,1\}$-valued symmetric matrix with zeros on the diagonal, it is the adjacency matrix of a simple graph, say $G$.  Our problem then reduces to finding the least positive integer $\lam$ such that the $\lam$-fold complete graph $\lam K_n$ admits a `signed' edge-decomposition into copies of $G$.  This question has been studied before \cite{CC,DL,Wilson75,WW}, although to our knowledge an explicit formula for $\lam$ has not been written down in all cases.  As observed in these references, especially \cite{CC,Wilson75}, this $\lambda$ is also the eventual period for multiples of the complete graph to admit an edge-decomposition (in the traditional `unsigned' setting) into copies of $G$.

The following special case of our work to follow settles a question of Cameron and Cioab\v{a} in \cite{CC}.  Independently, this can be derived from Wilson and Wong's Smith forms \cite{WW} of certain related inclusion matrices, although our proof is somewhat more direct.

\begin{thm}[see also \cite{WW}]
\label{thm:graphs}
Let $G$ be a simple graph with $n$ vertices, $e>0$ edges, and vertex degrees $d_1,\dots,d_n$.  
Put $\alpha=0$ if $G$ is regular; otherwise $\alpha=(2e-rn)/q$, where 
$q=\gcd(d_i-d_j: 1 \le i,j \le n)$ and $r$ is the least residue of each $d_i \pmod{q}$.
A signed $G$-decomposition of $\lam K_n$ exists if and only if $\lam_\text{min} \mid \lam$, where
\begin{equation}
\label{lambda-primitive}
\lam_\text{min} = \ell(G) := \frac{2e}{\gcd(n(n-1),\alpha(n-1),2e)},
\end{equation}
unless $G=K_{a,n-a}$ or $K_a \dotcup K_{n-a}$ for $n \neq 2a$, in which case 
\begin{equation}
\label{lambda-special}
\lam_\text{min} = \lcm(2,\ell(G)) = \frac{2e}{\gcd(n(n-1),a(n-1),a(n-a))}.
\end{equation}
\end{thm}

Our general problem on matrices can be recast into a similar setting, except that instead of $G$ being a simple graph, it is allowed to be a (possibly directed) multigraph equipped with integer edge weights, and possibly with vertex loops also having integer weights.  It is sometimes convenient to refer to terminology from graphs; for instance, a sum of off-diagonal entries in some row of our matrix may be thought of as a vertex outdegree.

Although we do not obtain a closed form as in Theorem~\ref{thm:graphs} for the most general version of our problem, the solution comes from a small integer lattice calculation.  We also give a simple constructive proof of existence of a matrix balancing when the necessary lattice condition holds.

The outline of the paper is as follows.  In the next section, we set up some notation and cover various preliminary results.  Among these is the analysis of general $3 \times 3$ matrices, which permits a useful assumption $n \ge 4$ in the sequel. In Section 3, we restrict attention to symmetric matrices with vanishing diagonal, which we model as signed multigraphs.  We prove a decomposition lemma which preprocesses a given signed multigraph into convenient summands.  This leads to a simple necessary and sufficient lattice condition for this version of our problem.  Theorem~\ref{thm:graphs} is obtained as a consequence.  Section 4 gives a similar treatment for general square matrices, where nonzero diagonal entries are allowed and the assumption of symmetry is dropped.  We conclude in Section 5 with a summary and brief look at a few next directions for the research.

\section{Preliminaries}
\label{sec:prelim}

\subsection{Notation and simple algebraic facts}

Let $\Z^{m \times n}$ denote the set of $m \times n$ integer matrices.  Let $\Sym_n$ denote the symmetric group on $[n]:=\{1,2,\dots,n\}$.

Given $A \in \Z^{n \times n}$, let $A^\sigma$ be the $n \times n$ matrix with entries
\begin{equation}
\label{Asigma}
(A^\sigma)_{ij} =A_{\sigma(i),\sigma(j)}.
\end{equation}
Equivalently, $A^\sigma = P^\top AP$, where $P$ is the permutation matrix for $\sigma$.
A matrix $A$ satisfies $A^\sigma = A$ for every $\sigma\in \Sym_n$ if and only if it is completely symmetric, i.e. $A \in \langle I,J \rangle := \{aI+bJ:a,b \in \Z\}$.

\begin{defn}
\label{defn:bal}
For $A \in \Z^{n \times n}$, the \emph{balancing index} or \emph{balancing number} of $A$, denoted $\bal(A)$, is the least positive value of 
$\sum_{\sigma \in \Sym_n} c_\sigma$,
where $(c_\sigma:\sigma \in \Sym_n)$ are integer coefficients such that
\begin{equation}
\label{Ccond}
\sum_{\sigma \in \Sym_n} c_\sigma A^\sigma \in \langle I,J \rangle.
\end{equation} 
\end{defn}
Observe that the set of possible values of $\sum_{\sigma \in \Sym_n} c_\sigma$ satisfying \eqref{Ccond} is a nontrivial ideal in $\Z$; so $\bal(A)$ is well-defined as the (positive) generator for that ideal.  We also remark that any completely symmetric combination of the form \eqref{Ccond} is a multiple of the matrix $\overline{A}=\frac{1}{n!} \sum_{\sigma \in \Sym_n} A^\sigma$, which has entries $\frac{1}{n} \mathrm{tr}(A)$ on the diagonal and 
$\frac{1}{n(n-1)} \mathrm{sum}^*(A)$ off the diagonal, where $\mathrm{sum}^*(A):=\sum_{i \neq j} A_{ij}$.

We consider some easy cases of the definition.  If $A$ is itself completely symmetric, then $\bal(A)=1$. However, as we see later the converse is not true in general.  For $2 \times 2$ matrices $A$, $\bal(A) =1$ or $2$, the former if and only if $A$ is completely symmetric.  If $B=P^\top AP$ is similar to $A$ via a permutation matrix, then $\bal(B)=\bal(A)$.
If $A$ is an $n \times n$ diagonal matrix, then $\bal(A) \mid n$, since the sum of copies $A^\sigma$ over the cyclic subgroup generated by $(1 2 \cdots n)$ satisfies \eqref{Ccond}.  This holds in various other cases, e.g. for back-circulant matrices of prime order.

Let $G$ be a graph with $V(G)=[n]$.
In what follows, it is convenient to slightly abuse notation, letting `$G$' also stand for its adjacency matrix.  In view of \eqref{Asigma}, it is natural to let $G^\sigma$ denote a copy of $G$ permuted by $\sigma \in \Sym_n$ on the same vertex set.  Our notation also facilitates the standard arithmetic on multigraphs by combining edge multiplicities.  Letting $G(\{i,j\})$ denote the multiplicity of edge $\{i,j\}$ in $G$, we have
\begin{enumerate}
\item
$G^\sigma(\{i,j\}) = G(\{\sigma(i),\sigma(j)\})$ for $V(G)=[n]$, $\sigma \in \Sym_n$;
\item
$(G_1+G_2)(\{i,j\})=G_1(\{i,j\})+G_2(\{i,j\})$ for $V(G_1)=V(G_2)=[n]$; and
\item
$(kG)(\{i,j\}) = k G(\{i,j\})$ for a number $k$ and multigraph $G$.
\end{enumerate}

To illustrate the notation, we have the relation 
\begin{equation}
\label{sn}
\sum_{\sigma \in \Sym_n} G^\sigma = 2(n-2)! |E(G)| K_n,
\end{equation}
since every edge of the complete graph $K_n$ appears equally often as an edge of some $G^\sigma$.

Let us also apply Definition~\ref{defn:bal} to a graph $G$, where again $G$ is used in place of its adjacency matrix.  We immediately notice some motivation for the `index' terminology, namely a relationship between $\bal(G)$ and the index of its automorphism group $\Aut(G)$ as a subgroup of $\Sym_n$.

\begin{lemma}
For a graph $G$ of order $n$, $\bal(G)$ divides $|\Sym_n:\Aut(G)|$.
\end{lemma}

\begin{proof}
Let $T$ be a set of right coset representatives for $\Aut(G)$ in $\Sym_n$.  
Working from \eqref{sn}, we have 
\begin{equation*}
2(n-2)! |E(G)| K_n = \sum_{\sigma \in \mathcal{S}_n} G^\sigma 
= \sum_{\alpha \in \Aut(G)} \sum_{\tau \in T}  G^{\alpha\tau} = |\Aut(G)| \sum_{\tau \in T} G^\tau.
\end{equation*}
Therefore, 
$\sum_{\tau \in T} G^\tau = \lambda K_n$, where $\lambda = 2(n-2)! |E(G)|/|\Aut(G)|$.  We note that $\lambda$ is an integer because the sum of all $G^\tau$ is an integer linear combination of graphs.  The result now follows from the definition, since
$\sum_{\tau \in T} 1 = |T|=|\Sym_n:\Aut(G)|$.
\end{proof}

Of course, the same result holds for multigraphs, directed graphs, graphs with loops, or general integer matrices under a suitably extended notion of automorphism group.

As can be expected, computation of $\bal(A)$ requires some light number theory.  For a prime $p$ and an integer $n$, let $\nu_p(n)$ denote the largest exponent $t \ge 0$ such that $p^t \mid n$, which we take to be $\infty$ if $n=0$.  A congruence $a \equiv b \pmod{k}$ is to be interpreted as equality $a=b$ when $k=0$.

\subsection{Smith normal form and integer linear systems}

Given $M \in \Z^{m \times n}$, there exist \cite{Hungerford,Newman} unimodular matrices $U \in \Z^{m \times m}$ and $V \in \Z^{n \times n}$ such that $UMV$ equals an $m \times n$ diagonal matrix, say $D$ with diagonal entries $d_1,d_2,\dots,d_{\min(m,n)}$, satisfying $d_i \mid d_{i+1}$ for each $i$.  The matrix $D$ is unique and called the \emph{Smith normal form} of $M$.
The surveys \cite{Newman,Sin} showcase some interesting combinatorial applications of Smith forms, e.g. to distinguish inequivalent Hadamard matrices or to spectral properties of strongly regular graphs.
Wilson and Wong \cite{WW} use the Smith form to study the module generated by copies of $G$ for 
simple graphs and certain special multigraphs $G$.

Given a system of linear equations $M \vx = \vb$, the Smith form of $M$ is useful to probe the existence of and compute integer solutions $\vx$.  Letting $\vy = V\vx$, an equivalent system is 
\begin{equation}
\label{equiv-system}
D \vy = U^{-1} \vb.
\end{equation}
A solution $\vy$ to \eqref{equiv-system} exists if and only if each $d_i$ divides the $i$th entry of $U^{-1}\vb$, in which case $y_i$ is the quotient.  Since $V$ is invertible over the integers, a solution $\vx$ can be recovered from $\vy$.

The relevant matrix $M$ for our problem of computing $\bal(A)$ is the $n^2 \times n!$ matrix whose columns are $\mathrm{vec}(A^\sigma)$, $\sigma \in \mathcal{S}_n$, where 
$\mathrm{vec}:\R^{n \times n} \rightarrow \R^{n^2}$ is the `vectorization' transformation, listing the elements of the input matrix as a column vector.  The ordering of the columns is unimportant.
This matrix is closely related to the matrix $N_2$ defined in \cite{WW}, which applies to graphs (and multigraphs).  In that setting, the matrix is a generalized inclusion matrix of edges of $K_n$ versus graph copies.  In view of \eqref{equiv-system}, having access to the Smith form of $N_2$ (or its matrix analog) greatly facilitates the computation of $\bal(A)$.  We let $\vb$ take the role of an unknown integer multiple of $\mathrm{vec}(\overline{A})$, at which point $\bal(A)$ is the least positive multiple allowing an integer solution to the resulting system \eqref{equiv-system}.

\subsection{Graph designs}

A graph $K$ is said to have a $G$-\emph{decomposition} if $E(K)$ can be partitioned into subgraphs  isomorphic to $G$.  The problem of determining sufficient conditions on $K$ to admit a $G$-decomposition has received considerable attention in recent years.  For instance, Delcourt and Postle show in \cite{DP} that sufficiently large graphs $K$ admit a $K_3$-decomposition provided the na\"ive arithmetic conditions hold, i.e. $3 \mid |E(K)|$ and $2 \mid \deg_K(x)$ for every $x \in V(K)$, and subject to minimum degree assumption $\delta(K) \ge 0.827|V(K)|$.

The special case when $K$ is a complete graph (or a multiple thereof) is more classical and originates with work of R.M.~Wilson, \cite{Wilson75}.  Motivated by a close connection to block designs, we say that a $G$-\emph{design} is a collection of copies of $G$ (called blocks) whose union is $\lam K_n$.  Note that we allow $|V(G)|<n$; isolated vertices can be included or ignored since they have no effect on edge decompositions. Wilson showed that, for fixed graphs $G$, sufficiently large complete graphs admit a $G$-decomposition whenever the global and local divisibility conditions hold.
An important step in the proof, and of primary interest for us here is the integral relaxation.  A \emph{signed} $G$-\emph{design} is an integer linear combination of copies of $G$ placed on vertices in $[n]$ that equals $\lam K_n$ under the edge arithmetic described earlier.  The balancing index $\bal(G)$ is effectively measuring the (net) number of blocks in a minimum nontrivial $G$-design.

\begin{prop}[Wilson, \cite{Wilson75}]
\label{prop:2isolates}
Let $G$ be a graph of order $n$ with $e$ edges, $e>0$, and at least two isolated vertices.
Let $d$ be the gcd of all degrees in $G$.
Put $$\lam_{\min} = \lcm\left( \frac{2e}{\gcd(n(n-1),2e)},\frac{d}{\gcd(n-1,d)} \right).$$
Then $\lam K_n$ admits a signed $G$-design if and only if $\lam_{\min} \mid \lam$.  That is,
$\bal(G) = \lam_{\min} n(n-1)/2e$.
\end{prop}

Working from Proposition~\ref{prop:2isolates}, Wilson observes that $\lam K_n$ actually has a $G$-decomposition (in the unsigned sense) for all sufficiently large multiples $\lam$ of $\lam_{\min}$.  To sketch the proof, $\lam_1 = 2(n-2)!|E(G)|$ admits a $G$-design via \eqref{sn}, as does $\lam_2 = \lam_{\min}+k \lam_1$ if we choose $-k$ as the smallest (negative) coefficient used in a signed $G$-design for $\lam_{\min}$.  Since $\gcd(\lam_1,\lam_2) = \lam_{\min}$, Frobenius' theorem delivers the conclusion.

It is natural to extend the notion of $G$-designs, both signed and unsigned, to the setting of square integer matrices.  We define a \emph{signed} $A$-\emph{design} as a list of integer coefficients $(c_\sigma)$ satisfying \eqref{Ccond}, and we drop the qualifier `signed' if the $c_\sigma$ are all nonnegative.
By exactly the same reasoning as above, $\bal(A)$ is not only the least positive coefficient total in a signed $A$-design, but also the eventual period for the size of large $A$-designs.

The simplifying assumption in Proposition~\ref{prop:2isolates} that $G$ have two isolated vertices is better understood in Section~\ref{sec:multigraph}.  The underlying reason is that signed $G$-designs can be built with the aid of gadgets of the form $G-G^\sigma-G^\tau+G^{\sigma\tau}$ for disjoint transpositions $\sigma,\tau \in \Sym_n$.  This leads us to consider the case $n=3$ separately, which we do next.

\subsection{The $3 \times 3$ case}
\label{sec:3x3}

Here, we restrict to $A \in \Z^{3 \times 3}$.  From our earlier remarks, it follows that $\bal(A) \in \{1,2,3,6\}$.  Deciding the value of $\bal(A)$ can be accomplished with the aid of the following two lemmas.

\begin{lemma}
\label{lem:parity3x3}
We have $2 \mid \bal(A)$ if and only if $A_{12}+A_{23}+A_{31} \neq A_{21}+A_{32}+A_{13}$.
\end{lemma}

\begin{proof}
Suppose first that $\bal(A)$ is odd.  Then there exists a $6$-tuple of coefficients $(c_\sigma)$ such that 
$B:=\sum_{\sigma \in \mathcal{S}_3} c_\sigma A^\sigma \in \langle I,J \rangle$ but where
\begin{equation}
\label{eqn:parity}
c_{()}+c_{(123)}+c_{(132)} \neq c_{(12)}+c_{(13)}+c_{(23)}.
\end{equation}
The sum of the $(1,2)$, $(2,3)$, and $(3,1)$-entries of $B$ is
$$(c_{()}+c_{(123)}+c_{(132)})(A_{12}+A_{23}+A_{31}) + 
(c_{(12)}+c_{(13)}+c_{(23)})(A_{21}+A_{32}+A_{13}).$$
The sum of the $(2,1)$, $(3,2)$, and $(1,3)$-entries of $B$ is
$$(c_{()}+c_{(123)}+c_{(132)})(A_{21}+A_{32}+A_{13})+ 
(c_{(12)}+c_{(13)}+c_{(23)})(A_{12}+A_{23}+A_{31}).$$
Equate these, since $B=B^\top$, and use \eqref{eqn:parity} to conclude 
$A_{12}+A_{23}+A_{31} = A_{21}+A_{32}+A_{13}$.

Conversely, if $A_{12}+A_{23}+A_{31} = A_{21}+A_{32}+A_{13}$, it is easy to see that 
$$A +A^{(123)}+A^{(132)} \in \langle I,J \rangle.$$
Thus $\bal(G) \mid 3$, and hence $\bal(G) \in \{1,3\}$.
\end{proof}

Testing whether $3 \mid \bal(A)$ is more subtle, but we are able to obtain a direct characterization.  We first consider the case of diagonal matrices.

\begin{lemma}
\label{ternary-triple}
The following are equivalent for integers $a,b,c$:
\vspace{-10pt}
\begin{enumerate}
\item[(a)]
$\bal(\mathrm{diag}(a,b,c)) =1$;
\item[(b)]
there exist integers $x,y,z$ such that $x+y+z=1$ and $ax+by+cz = \frac{1}{3}(a+b+c)$;
\item[(c)]
$\nu_3(a+b-2c) \ge 1+ \min(\nu_3(a-c), \nu_3(b-c))$.
\end{enumerate}
\end{lemma}

\begin{proof}
(a) $\Leftrightarrow$ (b): For the forward implication, consider any diagonal entry in an integer linear combination of $D^\sigma$ equaling $\overline{D}$.  Conversely, 
$$x(D+D^{(23)})+y(D^{(123)}+D^{(12)}) + z (D^{(132)}+D^{(13)}) = \begin{bmatrix}
u & 0 & 0\\
0 & v & 0\\
0 & 0 & v
\end{bmatrix},$$
where $u=2ax+2by+2cz = \frac{2}{3}(a+b+c)$ and $u+2v=2(x+y+z) \tr(D) =2(a+b+c)=3u$.
It follows that $\bal(D) \mid 2$, so $\bal(D) =1$ by Lemma~\ref{lem:parity3x3}.

(b) $\Leftrightarrow$ (c): 
Let $j=\min(\nu_3(a-c), \nu_3(b-c))$. Since $ax+by+cz = \frac{1}{3}(a+b+c)$, it follows that  $3(a-c)x+3(b-c)y+3(x+y+z)c=a+b+c$, or simply $3(a-c)x+3(b-c)y=a+b-2c$. By assumption $3^{j+1}$ divides the left-hand side proving the forward implication. 
%%% NEW 
Conversely, from $\gcd(a-c,b-c) \mid (a+b-2c)$ and (c) it follows that $3 \gcd(a-c,b-c) \mid (a+b-2c)$.
So there exist $x,y \in \Z$ such that $3x(a-c)+3y(b-c)=a+b-2c$.  Adding $3c$ to both sides, we obtain 
$3ax+3by+3c(1-x-y)=a+b+c$, as required.
%%% END NEW
\end{proof}

\rk
Condition (c) can be alternatively phrased as follows: for distinct $a,b,c$, they must all disagree in the the least significant ternary digit where they don't all agree.  We say that $(a,b,c)$ is a \emph{ternary triple} if it satisfies any of the above conditions.  The number of ternary triples $(a,b,c)$ with $1 \le a<b<c \le n$ is a new submission to the OEIS database \cite{OEIS} at \url{https://oeis.org/A349216}.

Let $d^+_i$, $d^-_j$ denote the $i$th row sum and $j$th column sum of $A$, respectively.
Define an auxiliary $4 \times 3$ matrix
$$\Theta(A)=
\left[
\begin{array}{ccc}
1 & 1 & 1 \\
A_{11} & A_{22}& A_{33}\\
d_1^+ & d_2^+ &d_3^+ \\
d_1^- & d_2^- &d_3^- 
\end{array}
\right].$$
Using $\Theta(A)$ and ternary triples, we can classify
those $3 \times 3$ matrices with balancing index $1$ or $2$.

\begin{prop}
\label{prop:3x3}
Let $A \in \Z^{3 \times 3}$.  We have $\bal(A) \le 2$ if and only if 
\vspace{-10pt}
\begin{enumerate}
\item[(1)]
$\mathrm{rank}(\Theta(A)) < 3$; and
\item[(2)]
the row sums, column sums, and diagonal entries of $A$ are ternary triples.
\end{enumerate}
\end{prop}

\begin{proof}
Suppose first that $\bal(A) \le 2$.  Let $d^+_i$, $d^-_j$ be the line sums of $A$, as above.  
Observe that the line sums of $A^\sigma$ are given by $d^+_{\sigma(i)}$ and $d^-_{\sigma(j)}$.  
Let $c_\sigma$ be coefficients with $\sum c_\sigma=2$ and $\sum c_\sigma A^\sigma = 2 \overline{A}$.  Reading this by rows, columns, and diagonal entries, we have
$\bal( \mathrm{diag}(d^+_i))$,  $\bal( \mathrm{diag}(d^-_j))$, and $\bal( \mathrm{diag}(A_{ii}))$ each at most $2$.
Lemma~\ref{lem:parity3x3} implies these numbers equal one, and hence condition (2) holds by 
Lemma~\ref{ternary-triple}.
Now consider the linear system
$$\Theta(A) \begin{bmatrix} x\\y\\z \end{bmatrix} = 
\begin{bmatrix} 2\\ \frac{2}{3} \tr(A)\\
\frac{2}{3} \mathrm{sum}(A) \\ 
\frac{2}{3} \mathrm{sum}(A)
\end{bmatrix}.$$
It has the solution $(x,y,z)=(\frac{2}{3},\frac{2}{3},\frac{2}{3})$ and also the integer solution
$(c_{()}+c_{(23)},c_{(123)}+c_{(12)},c_{(132)}+c_{(13)}).$ Their difference is a nonzero vector in the (right) kernel of $\Theta(A)$.

Conversely, suppose (1) and (2) hold.  We define a 
$\Z$-linear combination $A'$ of copies of $A$.
If rank$(\Theta(A))=1$, we let $A'=2A$. Otherwise, one of row $2$, $3$ or $4$ of $\Theta(A)$ is non-constant.  Suppose not all $A_{ii}$ are equal; the other cases are similar.   Then the column space of 
$\Theta(A)$ contains every ordered pair in the first two coordinates. From (1), a vector in the column space of $\Theta(A)$ is uniquely determined by its first two coordinates.  
Using (2), there exist integers $x,y,z$ such that $x+y+z=1$ and $A_{11} x + A_{22} y + A_{33} z =\frac{1}{3}\tr(A)$.  These agree with the first two coordinates of $\Theta(A)$ times $(\frac{1}{3},\frac{1}{3},\frac{1}{3})$.  So it follows that 
$$d_1^+x +d_2^+y + d_3^+z=d_1^-x +d_2^-y + d_3^-z=\tfrac{1}{3}\mathrm{sum}(A).$$
Put $c_{()}=c_{(23)}=x$, $c_{(123)}=c_{(12)}=y$, $c_{(132)}=c_{(13)}=z$, and let $A'=\sum_\sigma c_\sigma A^\sigma$. It is easy to check that $A'_{ii}=\frac{2}{3}\tr(A)$ for each $i=1,2,3$.  After a similar calculation, the first row sum of $A'$ is $2(d_1^+x +d_2^+y + d_3^+z)=\frac{2}{3}\mathrm{sum}(A)$ and both of the other row sums of $A'$ are $\mathrm{sum}(A)-(d_1^+x +d_2^+y + d_3^+z)=\frac{2}{3}\mathrm{sum}(A)$. So $A'$ has equal row sums.  Similarly, $A'$ has equal column sums.
From these properties, the off-diagonal entries satisfy
$$A'_{12}+A'_{13}=A'_{21}+A'_{23}=A'_{31}+A'_{32}=A'_{21}+A'_{31}=A'_{12}+A'_{32}=A'_{13}+A'_{23}.$$
This implies $A'_{12}=A'_{23}=A'_{31}$ and 
$A'_{21}=A'_{32}=A'_{13}$.  It follows that $A'=(A')^{(123)}$.  So $A'+(A')^{(12)} = 2\overline{A'}=4\overline{A}$. The sum of coefficients used on copies $A^\sigma$ in this combination equals $2(2x+2y+2z)=4$.  Therefore, $\bal(A) \mid 4$.  But $\bal(A) \mid n!=6$ and hence $\bal(A) \mid 2$.
\end{proof}

\section{The multigraph case}
\label{sec:multigraph}

In this section, we consider symmetric $n \times n$ matrices with a diagonal of zeros, and assume $n \ge 4$.  This includes the study of (undirected) multigraphs, although negative integer edge multiplicities are permitted with no extra difficulty.  Consider then a (signed) multigraph $G$ of order $n \ge 4$, say arising from the symmetric matrix $A$.  The `number of edges' of $G$, denoted $e$, is half of the sum of entries of $A$, and the degree $d_i$ of vertex $i$ is the $i$th row sum (or column sum) of $A$.  We call $(d_1,d_2,\dots,d_n)$ the \emph{degree vector} of $G$.

\subsection{Primitivity}

Let $M_G$ be the $\Z$-module generated by all permuted copies of $G$ on vertex set $[n]$. 
For an explicit presentation, as with matrices one can vectorize a signed multigraph into a standard presentation in $\Z^{\binom{n}{2}}$ by listing edge multiplicities (say) following a colex ordering.

Take disjoint transpositions, say $\sigma=(12)$, $\tau=(34)$ in $\Sym_n$, and consider 
$G^\Box := G - G^\sigma - G^\tau + G^{\sigma \tau} \in M_G$.
The multiplicity of $\{1,3\}$ in this linear combination is
\begin{equation}
\label{pod}
G^\Box(\{1,3\}) = G(\{1,3\})-G(\{2,3\})-G(\{1,4\})+G(\{2,4\}).\end{equation}
It is easy to see that the multiplicity of $\{2,4\}$ is also given by \eqref{pod}, of $\{2,3\}$ and $\{1,4\}$ is negative this value, and that all other multiplicities vanish.  In other words, $G^\Box$
equals \eqref{pod} times the `alternating $4$-cycle' $C_4^*$, whose vector form is given by
\begin{equation}
\label{c4}
\begin{tabular}{cccccccccc}
& \tiny $\{1,2\}$ 
& \tiny $\{1,3\}$ 
& \tiny $\{2,3\}$ 
& \tiny $\{1,4\}$ 
& \tiny $\{2,4\}$ 
& \tiny $\{3,4\}$ 
& \tiny $\{1,5\}$ 
\\
$($ & $0$, & $1$, & $-1$, & $-1$,& $1$, & $0$,& $0$,&  $\dots$ & $)$.
\end{tabular}
\end{equation}

Following \cite{WW}, the \emph{index of primitivity} of $G$ is 
$$k = \gcd( G^{\pi\Box}(\{1,3\}) : \pi \in \Sym_n),$$
and we say that $G$ is \emph{primitive} if $k=1$.  
By a straightforward argument, a simple graph $G$ of order $n \ge 4$ is primitive unless it is isomorphic to $K_{a,b}$ or its complement $K_a \dotcup K_b$. Here, $(a,b) \in \{(0,n),(n,0)\}$ is understood to reference an edgeless or complete graph.  The assumption of 2 isolates in Proposition~\ref{prop:2isolates} ensures primitivity for simple graphs, though this condition is in general far stronger than needed.  Details on this case can be found in \cite{Wilson75,WW}.

Note that it is possible to have $k=0$ in nontrivial cases.  A combination $G^\Box$ is also known as a $2$-pod; these are discussed in more detail and generality in \cite{GJ,WW}.

Before moving on to general multigraphs, we briefly examine the primitive ones.  Here, $\bal(G)$ can in principle be computed from the Smith forms in \cite{WW}, but we present an alternate approach for which the primitive case is especially easy.  Suppose $G$ is primitive with degree vector $\vd$.  A signed $G$-decomposition of $\lam K_n$ can exist only if the submodule of $\Z^n$ generated by all permutations of $\vd$ contains  $\lam(n-1) \one$, the degree vector of $\lam K_n$.  Conversely, supposing this lattice condition holds, one can begin with any degree-balanced combination of copies of $G$ and iteratively reduce the edge variance using copies of $C_4^*$ until the resulting multigraph has zero edge variance, i.e. is of the form $\lam K_n$.  More details are given later, but this serves to outline the method and motivates building (multiples of) $C_4^*$ in $M_G$.

\subsection{Coboundaries}
\label{sec:coboundaries}
To handle general index of primitivity, we define a `coboundary operator' $\delta$ as follows.
Given a function $f:[n] \rightarrow \R$, we have $\delta f : \binom{[n]}{2} \rightarrow \R$ defined by
\begin{equation*}
(\delta f)(\{i,j\}) = f(i) + f(j).
\end{equation*}
It is easy to see that $\delta$ is a linear transformation, and that $\text{range}(\delta f) \subset \Z$, i.e. $\delta f$ is a signed multigraph, if and only if $\text{range}(f)$ is contained either in $\Z$ or $\frac{1}{2} \Z \setminus \Z$.  A \emph{coboundary graph} is a signed multigraph in the image of $\delta$.  That is, a coboundary graph is determined by an assignment of \emph{vertex potentials} $f_1,\dots,f_n$, either all in $\Z$ or all in $\frac{1}{2} \Z \setminus \Z$, and the multiplicity of edge $\{i,j\}$ equals $f_i+f_j$.  Given a coboundary graph $F$ with potentials $f_1,\dots,f_n$, its `complement' $K_n-F$ is another coboundary graph with potentials given by $f_i'=\frac{1}{2}-f_i$ for each $i$.

\begin{figure}[htbp]
\begin{center}
\begin{tikzpicture}
\draw [red,dashed] (0:1)--(72:1)--(288:1)--(0:1);
\draw (144:1.05)--(216:1.05);
\draw (144:0.95)--(216:0.95);
\draw (216:1)--(72:1);
\draw (216:1)--(0:1);
\draw (216:1)--(-72:1);
\node at (0:1.35) {\small $-\frac{1}{2}$};
\node at (72:1.3) {\small $-\frac{1}{2}$};
\node at (-72:1.3) {\small $-\frac{1}{2}$};
\node at (144:1.3) {\small $\frac{1}{2}$};
\node at (216:1.3) {\small $\frac{3}{2}$};
\foreach \a in {0,1,...,4}
 \filldraw (\a*72:1) circle [radius=0.1];
\end{tikzpicture}
\hspace{1cm}
\begin{tikzpicture}
\draw (0:1.05)--(72:1.05)--(288:1.05)--(0:1.05);
\draw (-72:0.95)--(0:0.95)--(72:0.95);
\draw (-74:0.9)--(74:0.9);
\draw [red,dashed] (144:1.05)--(216:1.05);
\draw (144:1)--(72:1);
\draw (144:1)--(0:1);
\draw (144:1)--(-72:1);
\node at (0:1.3) {\small $1$};
\node at (72:1.3) {\small $1$};
\node at (-72:1.3) {\small $1$};
\node at (144:1.3) {\small $0$};
\node at (216:1.35) {\small $-1$};
\foreach \a in {0,1,...,4}
 \filldraw (\a*72:1) circle [radius=0.1];
\end{tikzpicture}
\caption{A coboundary graph and its complement labelled with vertex potentials}
\label{fig:coboundary}
\end{center}
\end{figure}
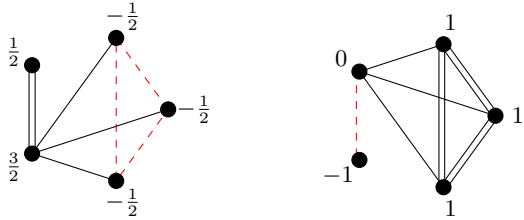

\begin{ex}
In Figure~\ref{fig:coboundary}, a coboundary graph and its complement are shown, with vertex potentials as indicated.  Dashed edges represent a multiplicity of $-1$.
\end{ex}

As mentioned earlier, nonprimitive graphs $G$ present a challenge that $C_4^* \not\in M_G$.  An illustrative case is $G=K_a \dotcup K_b$, a disjoint union of two cliques, say on vertex sets $A$ and $B$.  When $a,b \ge 2$, it is easy to see that this graph has index of primitivity equal to $2$.  But notice that if we subtract from $G$ the coboundary graph $\delta (\frac12 \one_A-\frac12 \one_B)$, the resulting graph is $0K_a \dotcup 2K_b$, which is twice a primitive graph.  In general, we would like to decompose an arbitrary multigraph in such a way so that the summands are simpler to handle.  This is the purpose of the following result.

\begin{lemma}[Primitive decomposition]  
\label{lem:prim-dec}
Let $G$ be a signed multigraph with index of primitivity $k$. Then there exists a coboundary graph $F$ and a primitive graph $H$ of the same order such that $G=F+kH$.
\end{lemma}

\begin{proof}
Put $a=G(\{1,2\})+G(\{1,3\})-G(\{2,3\})$.  For $i \in [n]$, define vertex potentials
$$f_i := 
\begin{cases}
a/2 & \text{if }i=1,\\
G(\{1,i\}) - a/2 & \text{otherwise}.
\end{cases}$$
Let $F$ be the coboundary $\delta f$, and put $H':=G-F$.  
Consider two vertices $i<j$.  If $i=1$,  then
$$H'(\{i,j\}) = G(\{i,j\})-F(\{i,j\}) = G(\{1,j\})-G(\{1,j\})+a/2-a/2 = 0.$$
For $i>1$,
\begin{equation}
\label{ij1a}
H'(\{i,j\}) = G(\{i,j\})-F(\{i,j\}) = G(\{i,j\})-G(\{1,i\})-G(\{1,j\})+a.
\end{equation}
By definition of $a$, \eqref{ij1a} vanishes for $(i,j)=(2,3)$. 
When $i=2$ and $j>3$, \eqref{ij1a} becomes 
$$G(\{2,j\})-G(\{1,j\})+G(\{1,3\})-G(\{2,3\}) \equiv 0 \pmod{k},$$
where we have used transpositions $\sigma=(12)$, $\tau=(3j)$ and that $G$ has index of primitivity equal to $k$.  Finally, for $i>2$, we similarly have 
\begin{align*} 
G(\{i,j\})+G(\{1,2\})-G(\{1,j\})-G(\{i,2\}) &\equiv 0\pmod{k},\\
G(\{i,j\})+G(\{1,3\})-G(\{1,i\})-G(\{j,3\}) &\equiv 0\pmod{k},\\ 
G(\{i,j\})+G(\{2,3\})-G(\{i,2\})-G(\{j,3\}) &\equiv 0\pmod{k}.
\end{align*} 
Adding the first two congruences and subtracting the third produces the right side of \eqref{ij1a}; hence
$H'(\{i,j\}) \equiv 0 \pmod{k}$ in all cases.  If $k=0$, this shows $G=F$; otherwise, put $H = \frac{1}{k}H'$.  It remains to show that $H$ is primitive. 
Since $G$ has index of primitivity $k$,
$$\gcd(\{r:G^{\pi \Box}=rC_4^*~\text{for some }\pi \in \mathcal{S}_n\}).$$
Since $F$ is a coboundary, $F^{\pi\Box}=0$ for all $\pi \in \mathcal{S}_n$.  It follows that $H'=G-F$ also has index of primitivity $k$, and thus $H=\frac{1}{k}H'$ is primitive.
\end{proof}

\subsection{The local lattice}

Let $G$ be a multigraph on vertex set $[n]$ with `degree vector', i.e. list of vertex degrees $\vd=(d_1,d_2,\dots,d_n)$.  For $\sigma \in \Sym_n$, the list of degrees for $G^\sigma$ is 
$\vd^\sigma:=(d_{\sigma(i)}:i=1,2,\dots,n)$.
From Lemma~\ref{lem:prim-dec}, we can take a list of vertex potentials $\vf=(f_1,f_2,\dots,f_n)$ that induce a coboundary $F$ agreeing with $G \pmod{k}$.   Potentials for $G^\sigma$ are 
$\vf^\sigma :=(f_{\sigma(i)}:i=1,2,\dots,n)$.

Consider a signed $G$-decomposition of $\lam K_n$, say
\begin{equation}
\label{Gdec-nec}
\sum_{\sigma \in \Sym_n} c_\sigma G^\sigma = \lam K_n.
\end{equation}
Taking the degree vectors on both sides gives
\begin{equation}
\label{deg-nec}
\sum_{\sigma \in \Sym_n} c_\sigma \vd^\sigma = \lam (n-1) \one.
\end{equation}
Next, we apply Lemma~\ref{lem:prim-dec} to the left side of \eqref{Gdec-nec} and get
$$\sum_{\sigma \in \Sym_n} c_\sigma F^\sigma \equiv \lam K_n \pmod{k},$$
by which we mean that every edge multiplicity on the left differs from $\lam$ by a multiple of $k$.
The potentials in $\lam K_n$ are each $\lam/2 \pmod{k}$, or $\lam/2+k/2 \pmod{k}$ if $k$ is even.
It is convenient to use twice the potentials, giving
\begin{equation}
\label{pot-nec}
\sum_{\sigma \in \Sym_n} c_\sigma (2\vf^\sigma) \equiv \lam \one \pmod{k}.
\end{equation}

For a fixed $G$, combining \eqref{deg-nec} and \eqref{pot-nec} results in a divisibility lower bound on $\lam$ necessary for the existence of a signed $G$-decomposition of $\lam K_n$.  The relevant module
here is $\langle \mathbf{d}^\sigma \oplus \mathbf{f}^\sigma : \sigma \in \mathcal{S}_n \rangle$, which occupies $2n$ dimensions and has $n!$ generators. 
Fortunately, we can exploit symmetry to greatly simplify the analysis.

Define the \emph{local lattice} of $G$ as 
\begin{equation}
\label{local-lattice}
\Lambda^*_G:= \langle \{ (1,d_i,2f_i) : i =1,\dots,n\} \cup \{(0,0,k) \} \rangle 
\subseteq \Z^3.
\end{equation}

\comment{
In what follows, we use 
$\oplus$ to emphasize a division into degree coordinates and potential coordinates.  We also use $\mathbf{u}_j$ to denote the $j$th standard basis vector in $\R^n$.  
With this in place, define the $\Z$-module
\begin{equation}
\label{intermediate-lattice}
\Lambda_G:=\langle \{ \vd^\sigma \oplus \vf^\sigma : \sigma \in \Sym_n \} \cup 
\{ \mathbf{0} \oplus k \mathbf{u}_j : j=1,\dots,n\} \cup \{\mathbf{0} \oplus \tfrac{k}{\gcd(k,2)} \one \} \rangle \subseteq \Z^{2n}.
\end{equation}
The following result summarizes the above necessary conditions on signed $G$-designs.

\begin{prop}
\label{prop:G-intermediate}
If there exists a signed $G$-decomposition of $\lam K_n$, then
$\lam ( (n-1) \one \oplus  \frac{1}{2}\one) \in \Lambda_G$, where $\Lambda_G$ is as defined in \eqref{intermediate-lattice}.
\end{prop}

This is effectively a (divisibility) lower bound on $\bal(G)$, but as stated still requires a large lattice or Smith form calculation.  However, w
}

Strictly speaking, a lattice in $\Z^3$ must fill three dimensions, and it can happen in special cases that $\Lambda^*_G$ occupies only a proper subspace of $\R^3$.  In any event, we show in what follows that 
$\bal(G)$ can be computed using only $\Lambda^*_G$. 
First, we give a lower bound on $\bal(G)$.

%%% new stuff %%%%%%%%%%%%%%%%%%%%%%%%%%%%%%%%

\begin{prop}
\label{prop:G-local}
Let $G$ be a multigraph of order $n$ with $e$ edges and local lattice $\Lambda^*_G$.  If there exists a signed $G$-decomposition of $\lam K_n$ with $C$ blocks (the sum of coefficients in the decomposition), then
$(C, \lam(n-1), \lam) = C (1,2e/n,2e/n(n-1)) \in \Lambda^*_G$. The least positive such $C$ divides
$\bal(G)$.
\end{prop}

\begin{proof}
Suppose $c_\sigma$
are integer coefficients satisfying 
$\sum_{\sigma \in \mathcal{S}_n} c_\sigma = C$ and
$\sum_{\sigma \in \mathcal{S}_n} c_\sigma G^\sigma = \lam K_n$,
where $\lam=\frac{2eC}{n(n-1)}$. 
For $i \in [n]$, let $C_i$ be the sum of $c_\sigma$ over all $\sigma \in \mathcal{S}_n$ satisfying $\sigma(1)=i$. Then
$$\sum_{i=1}^n C_i = \sum_{\sigma \in \mathcal{S}_n} c_\sigma = C.$$
Next, examining degrees and potentials at vertex $1$,  \eqref{deg-nec} and \eqref{pot-nec} give
$$\sum_{i=1}^n C_i d_i = \sum_{\sigma \in \mathcal{S}_n} c_\sigma (\mathbf{d}^\sigma)_1 = \lam(n-1)$$
and
$$\sum_{i=1}^n C_i (2f_i) = \sum_{\sigma \in \mathcal{S}_n} c_\sigma (2\mathbf{f}^\sigma)_1 \equiv \lam \pmod{k}.$$
It follows that 
$$(C, \lam(n-1), \lam) = C (1,2e/n,2e/n(n-1))
= \sum_{i=1}^n C_i (1,d_i,2f_i) \in \Lambda^*_G.$$
The last claim follows by taking $C=\bal(G)$ in the above calculations.
\end{proof}

\subsection{Sufficiency}

Here, we prove that the above local lattice condition is sufficient, thereby determining $\bal(G)$ for general (signed) multigraphs.

To assist in this, we define for each $r \in [n]$ the group ring element
$$\alpha_r:=\sum_{s=1}^n (1sr) - (n-1) (1r) \in \Z[S_n].$$
In the special case $r=1$, this reduces to 
$\sum_{s=1}^n (1s) - (n-1) ()$.  To motivate this definition, consider a vector $\mathbf{v}=(v_1,\dots,v_n)$.
Under the natural action on coordinates, 
$$\mathbf{v}^{\alpha_r}=\left(\mathrm{sum}(v)-(n-1)v_r,v_r,\dots,v_r\right).$$
Loosely speaking, $\alpha_r$ isolates component $r$, and does so with a coefficient sum equal to $1$. This allows us to lift $\Z$-linear combinations in the local lattice $\Lambda^*_G$ to combinations in the full module $M_G$.

\begin{thm}
\label{thm:multigraphs}
Let $G$ be a multigraph of order $n$ with local lattice $\Lambda^*_G$, and let $C$ denote the least positive multiple of $(1,2e/n,2e/n(n-1))$ belonging to $\Lambda^*_G$.  Then $\bal(G)=C$.
\end{thm}

\begin{proof}
Proposition~\ref{prop:G-local} shows that $C \mid \bal(G)$, so it suffices to prove the converse.
Suppose $C_i$ are integer coefficients satisfying
\begin{equation}
\label{sufficient-comb}
\sum_{i=1}^n C_i(1,d_i,2f_i)+C_0(0,0,k) = C(1,2e/n,2e/n(n-1))=(C,\lam(n-1),\lam).
\end{equation}
Define the multigraph $K:=\sum_{i=1}^n C_i G^{\alpha_i} \in M_G$. 
Since each $\alpha_i$ has a coefficient sum of $1$, the total sum of coefficients in $K$ is $\sum_{i=1}^n C_i=C$. We claim that $K$ has two properties: (1) it is regular of degree $\lam(n-1)$; and (2) $K \equiv \lam K_n \pmod{k}$.

(1) Consider a vertex $x \in [n]$. The degree of $x$ in $G^{\alpha_i}$ equals $d_i$ if $x>1$ and is 
$2e-(n-1)d_i$ if $x=1$.  So, from \eqref{sufficient-comb}, $\deg_K(x)=\sum_{i=1}^n C_i d_i = \lambda(n-1)$ when $x>1$ and 
$\deg_K(1)=2eC-\lambda(n-1)^2=\lambda(n-1)$ as well.\\
(2) For distinct vertices $x$ and $y$, consider the edge multiplicity $K(\{x,y\})$.  
The primitive decomposition $G\equiv F \pmod{k}$ implies $K(\{x,y\}) \equiv \sum_{i=1}^n C_i F^{\alpha_i} (\{x,y\}) \pmod{k}$. When $x,y>1$, \eqref{sufficient-comb} gives
$$K(\{x,y\}) \equiv \sum_{i=1}^n C_i(f_i+f_i) \equiv \lam \pmod{k}.$$
To handle the case $x=1$, we can consider the sum over all $x\neq y$ to obtain
$$\lam(n-1) = \deg_K(y)=K(\{1,y\})+\sum_{x \neq 1,y} K(\{x,y\}) \equiv K(\{1,y\})+(n-2)\lam \pmod{k}.$$
It follows that $K(\{1,y\}) \equiv \lam \pmod{k}$ as well.

Now, among all multigraphs in $M_G$ which satisfy (1) and (2), assume $K$ minimizes the variance of edge multiplicities, 
$\text{var}(K):=\sum_{x<y} |K(\{x,y\})-\lam |$.
Suppose we have some pair $\{x,y\}$ such that $K(\{x,y\})>\lambda$. Then since 
$\deg_K(y)=\lambda(n-1)$, there is some $z\neq x$ such that $K(\{y,z\})<\lambda$. Similarly, since $\deg_K(z)=\lambda(n-1)$ there is some $w$ such that $K(\{z,w\})>\lambda$. If $w=x$, then we must have some $u \neq y$ such that $K(\{u,x\})<\lambda$.  Either way, we end up with a path of length three such that the edge multiplicities alternate between being strictly less than and greater than $\lambda$.  Let us relabel the vertices (in order) as $u,v,w,x$.

We know that $M_G$ contains a copy of $kC_4^*$ on vertices $u,v,w,x$.  Depending on signs, one of $K \pm kC^*_4$, has, compared with $K$, a decreased contribution to edge variance on each of $\{u,v\}$, $\{v,w\}$, and $\{w,x\}$ by $k$, and at worst increased the contribution on $\{x,u\}$ by $k$.  So, choosing a sign without loss of generality, 
$\text{var}(K + kC^*_4) < \text{var}(K)$.
Since this adjusted graph also satsifies (1) and (2), we arrive at a contradiction to the choice of $K$.
\end{proof}

\rk
The proof converts naturally into an algorithm that produces a signed $G$-decomposition from the local lattice condition. First, build a graph $K \in M_G$ using the group ring element $\sum C_i\alpha_i$.  Then, iteratively reduce variance on three-edge paths in $K$ until the resulting graph is $\lam K_n$.

We conclude with a summary restatement.

\begin{cor}
\label{lambda}
The following are equivalent for multigraphs $G$ of order $n \ge 4$.
\vspace{-10pt}
\begin{enumerate}
\item[(a)]
there exists a signed $G$-decomposition of $\lam K_n$;
\item[(b)]
$2e \, \bal(G)/n(n-1)$ divides $\lam$; and
\item[(c)]
$C(1,2e/n,2e/n(n-1))\in\Lambda^*_G$, where $\lam=2eC/n(n-1)$. 
\end{enumerate}
\end{cor}

\subsection{Some noteworthy cases}
\label{sec:noteworthy-graphs}

In the case when $G$ is a primitive multigraph, we can drop the 
`potential' coordinate from its local lattice $\Lambda^*_G$.  From Corollary~\ref{lambda}, the condition on $\lam$ reduces to $\lam (n(n-1),n-1) \in \langle (2e,d_i): i=1,\dots,n \rangle$.
With the help of a couple of easy lemmas, we can obtain a closed form in this case.

\begin{lemma}
\label{lem:2d-basis}
Let $d_1,\dots,d_n,m \in \Z$ with $m \neq 0$.  A basis for the $\Z$-module $\langle (m,d_i): i=1,\dots,n \rangle$ is $\{(m,d)\}$ if each $d_i=d$, and otherwise $\{(m,r),(0,q)\}$, where $q=\gcd(d_i-d_j: 1\le i,j\le n)$ and $d_i \equiv r \pmod{q}$ for each $i$.  
\end{lemma}

\rks
In the nondegenerate case, it is possible to take $r=d_1$, though if $r$ is chosen with $|r| \le q/2$ then one can obtain a minimal lattice basis.
Letting $d:=\gcd(d_1,\dots,d_n)$, we note that $d$ divides $q$.  For $G$ with an isolated vertex, (and for many generic graphs), one has $q=d$ and thus $r=0$.  In this situation, our basis is diagonal and the computation of $\bal(G)$ or equivalently $\lam_{\min}$ reduces to a straightforward least common multiple.  
On the other hand, it can happen that $d$ properly divides $q$ in general.  This simplification $q=d$ was mistakenly made in \cite[Lemma 2.1]{DL}, although elsewhere in that paper the proper lattice condition is used.

\begin{lemma}
\label{lem:2d-closedform}
For integers $m,q,r,u,v$, the least positive multiple $\ell$ of $(u,v)$ belonging to $\langle (m,r), (0,q) \rangle$ is 
$\ell=m/\gcd(m,u)$ if $q=0$ (provided $mv=ru$), and otherwise
\begin{equation*}
\ell =  \frac{mq}{\gcd(mq,uq,mv-ru)}.
\end{equation*}
\end{lemma}

Setting $m=2e$, $u=n(n-1)$, $v=n-1$ and cancelling a common factor of $q$, we obtain \eqref{lambda-primitive} as an expression for 
$\lam_{\min}$ in the primitive multigraph case.  Equivalently, using Theorem~\ref{thm:multigraphs}, 
\begin{equation}
\label{balG-primitive}
\bal(G) = \frac{n(n-1)}{\gcd(2e,n(n-1),\frac{2e-rn}{q}(n-1))},
\end{equation}
where the last term in the gcd is omitted if $G$ is regular.  This settles Theorem~\ref{thm:graphs} for primitive $G$.

As observed in \cite{WW}, the only simple graphs which are not primitive are of the form $K_a \dotcup K_b$ or $K_{a,b}$, where $a+b=n$ and $a \ge b \ge 0$. 
Consider the graph $G=K_{a,b}$, say with vertex partition $A \dotcup B$.  If $a \ge b \ge 2$, we have 
index of primitivity $k=2$ and primitive decomposition
$G= (0K_a \vee 2K_b) - 2K_b$ arising from potential vector $\vf=\one_B$.  The local lattice 
is
$\Lambda^*_G=\langle (1,b,0),(1,a,2),(0,0,2) \rangle$.  
If $b=1$, we simply drop the third vector, since then the index of primitivity equals zero.  
Ignoring potentials for a moment, a substitution $m=2ab$, $q=a-b$, $r=b$, $u=n(n-1)$, $v=n-1$
revovers $\ell(G)$ as in \eqref{lambda-primitive}, where we note that $\alpha = 0$ if $a=b$ and otherwise $\alpha=(2ab-bn)/(a-b) = b$.
But the even potentials imply that $\lam_{\min} = \lcm(\ell(G),2)$.  That is, balancing the potentials as well as the degrees forces an extra parity condition that $\lam$ be even. 

\begin{figure}[htbp]
\begin{center}
\includegraphics[width=5cm]{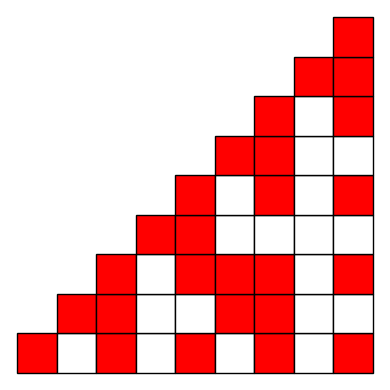}
\includegraphics[width=5cm]{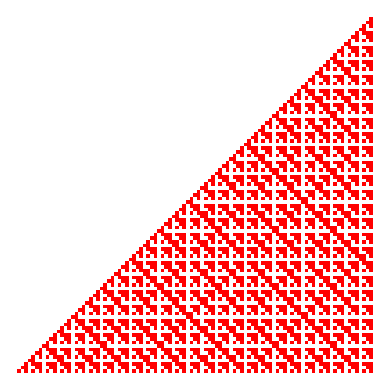}
\caption{Level set $\psi^{-1}(2)$ for $1 \le b \le a <10$ (left) and $<100$ (right).}
\label{fig:fractal}
\end{center}
\end{figure}

To study the parity in more detail, define 
$$\psi(a,b)=1+\left( \gcd\left( \frac{\lcm(a,b)}{\gcd(a,b)},
\frac{\lcm(a(a-1),b(b-1))}{\gcd(a(a-1),b(b-1))}  \right) \text{mod } 2 \right);$$
that is, $\psi(a,b)$ equals $2$ if $\nu_2(a)=\nu_2(b)$ or $\nu_2(a(a-1))=\nu_2(b(b-1))$, and otherwise $\psi(a,b)=1$. (If $b=1$, we ignore the latter $2$-adic test and obtain $\psi(a,1) =1+a \text{ mod } 2$.)   This function is shown in Figure~\ref{fig:fractal}, where cell $(a,b)$ is shaded if and only if $\psi(a,b)=2$; here $a=1,2,\dots$ indexes the horizontal axis and $b=1,2,\dots,a$ vertically.

Some easy case checking shows that $\ell(G) = 2ab/\gcd(n(n-1),a(n-1),2ab)$ is odd
if and only if $\psi(a,b)=2$.  It follows that $\lam_\text{min} =\psi(a,b) \ell(G)$ and
$\bal(G) = \psi(a,b)$ times the expression \eqref{balG-primitive}.

The other case $G=K_a \dotcup K_b$ can be studied similarly, or by complementing.  
Graphs $K_n$ and $\overline K_n$ are trivial, with $\lam=1$ and $0$, respectively.  This completes the analysis of simple graphs.  See also \cite[Section 7]{WW}.

Finally, we remark that if $G$ is a (signed) multigraph adorned with loops, say $l_i$ loops at vertex $i$, for $i=1,\dots,n$, then $\bal(G)$ can be computed by appending this extra data to the local lattice.  With this in mind, put
$$\Lambda_G^* =\langle \{(1,l_i,d_i,2f_i) : i = 1,\dots,n\} \cup \{(0,0,0,k)\} \rangle.$$
The `test vector' becomes $(1,p/n,2e/n,2e/n(n-1))$, where $p=l_1 + \dots + l_n$.
In the special case when $G$ has two isolated vertices (but possibly loops at other vertices), one can drop the first and third coordinates; this situation is covered in \cite[Lemma 2.3]{DM}.
In general, allowing loops essentially places us in the setting of general symmetric matrices $A \in \Z^{n \times n}$, where loop counts take the role of diagonal elements $A_{ii}$.
We cover this extension in more detail in the next section.

\section{General matrices}

Here, we extend results of the previous section from signed multigraphs to arbitrary (not necessarily symmetric) matrices $A \in \Z^{n \times n}$.  Again we assume $n \ge 4$, since smaller dimensions were discussed in Section~\ref{sec:prelim}.

\subsection{Special combinations and the primitivity lattice}
\label{sec:prim-lattice}

For symmetric matrices, the index of primitivity gives a concise description of the possible off-diagonal adjustments we can make of the form $\sum_\sigma c_\sigma A^\sigma$ where $\sum_\sigma c_\sigma = 0$.  Here, we expand this discussion to handle general matrices. 

Let $\mathcal{K}_4=\{ \mathrm{id},(12),(34),(12)(34) \} \le \Sym_n$ be the Klein $4$-group.  For $A \in \Z^{n \times n}$,
define 
$$A^{\Box} = \sum_{\sigma \in \mathcal{K}_4} \mathrm{sgn}(\sigma) A^{\sigma} = A-A^{(12)}-A^{(34)}+A^{(12)(34)}.$$
It is easy to see that
\begin{equation}
\label{K4-comb}
A^\Box = 
\left[
\begin{array}{rrrrrr}
0 & 0 & b & -b & 0 & \cdots \\
0 & 0 & -b & b & 0 & \cdots \\
c & -c & 0 & 0 & 0 & \cdots \\
-c & c & 0 & 0 & 0 & \cdots \\
0 & 0 & 0 & 0 & 0 & \\
\vdots & \vdots & \vdots & \vdots & & \ddots
\end{array}
\right]
\end{equation}
for some integers $b$ and $c$.  Consider the submodule $\Phi(A)$ of $\Z^2$ generated by all such pairs $(b,c)$ realized in $A^{\pi \Box}$, $\pi \in \mathcal{S}_n$.  We call $\Phi(A)$ the \emph{primitivity lattice} of $A$, since this generalizes the index of primitivity for symmetric $A$.  It is straightforward to extend Lemma~\ref{lem:prim-dec} to decompose general matrices $A$.

\begin{lemma}[Primitive decomposition]
\label{lem:matrix-prim-dec}
Given $A \in \Z^{n \times n}$, we can write
$A=B+F+Z$, where
$B$ is a matrix whose symmetric off-diagonal pairs $(B_{ij},B_{ji})$ all belong to $\Phi(A)$,
$F=[f_i^+ + f_j^-]_{ij}$ for some integers $f_i^+$ and $f_j^-$, and
$Z$ is a matrix which is constant above and constant below the main diagonal.
\end{lemma}

\begin{proof}
Put $z^+=A_{12}+A_{31}-A_{32}$, $z^-=A_{21}+A_{13}-A_{23}$, and for $i,j \in [n]$
\begin{equation}
\label{inout-pots}
f_i^+ := 
\begin{cases}
z^+ & \text{if }i=1,\\
A_{i1}  & \text{otherwise},
\end{cases}~\text{ and }
f_j^- := 
\begin{cases}
z^- & \text{if }j=1,\\
A_{1j} & \text{otherwise}.
\end{cases}
\end{equation}
Define $F,Z \in \Z^{n \times n}$ by
$F_{ij}=[f_i^+ + f_j^-]$ and 
$$Z_{ij}=
\begin{cases}
-z^+ & \text{if }i<j,\\
-z^- & \text{if }i>j,\\
0 & \text{otherwise}.
\end{cases}$$
Put $B=A-F-Z$.
We have $B_{i1} = B_{1j}= 0$ by \eqref{inout-pots}.
For $1<i<j$,
\begin{equation}
\label{ij1a-matrix}
B_{ij} = A_{ij}-A_{i1}-A_{1j}+z^+ \text{ and } B_{ji} = A_{ji}-A_{j1}-A_{1i}+z^-.
\end{equation}
For $i,j=(2,3)$, this gives $B_{23}=-B_{32}=A_{12}+A_{23}+A_{31}-A_{21}-A_{32}-A_{13}$.
Now, we have
\begin{align*}
(A_{12}-A_{32}-A_{14}+A_{34},A_{21}-A_{23}-A_{41}+A_{43}) & \in \Phi(A),\\
(A_{31}-A_{21}-A_{34}+A_{24},A_{13}-A_{12}-A_{43}+A_{42}) & \in \Phi(A),\\
(A_{23}-A_{13}-A_{24}+A_{14},A_{32}-A_{31}-A_{42}+A_{41}) & \in \Phi(A).
\end{align*}
Summing these, we find that $(B_{23},B_{32}) \in \Phi(A)$.  Similar relations imply
$(B_{ij},B_{ji}) \in \Phi(A)$ for any $i < j$.  Finally, $\Phi(A)$ is invariant under coordinate swap, so the same holds for $i>j$.
%\begin{align*}
%(A_{ij}+A_{12}-A_{i2}-A_{1j},A_{ji}+A_{21}-A_{2i}-A_{j1}) & \in \Phi(A),\\
%(A_{ij}+A_{31}-A_{i1}-A_{3j},A_{ji}+A_{13}-A_{1i}-A_{j3}) & \in \Phi(A),\\
%(A_{ij}+A_{32}-A_{i2}-A_{3j},A_{ji}+A_{23}-A_{2i}-A_{j3}) & \in \Phi(A).
%\end{align*}
\end{proof}
Note in particular that the proof shows
\begin{equation}
\label{z-phi}
(z^+-z^-)(1,-1)=(A_{12}+A_{23}+A_{31}-A_{21}-A_{32}-A_{13})(1,-1) \in \Phi(A).
\end{equation}
That is, $(z^+,z^-) \equiv (z^-,z^+) \pmod{\Phi(A)}$.
We use these facts later in the coming subsections.

\subsection{Triangles and a parity condition}
\label{sec:parity}

For $A \in \Z^{n \times n}$, define
\begin{equation}
\label{S3-comb}
A^{\triangle} = \sum_{\sigma \in \mathcal{S}_3} \mathrm{sgn}(\sigma) A^{\sigma} = 
\left[
\begin{array}{rrrrr}
0 & a & -a & 0 &  \cdots \\
-a & 0 & a & 0 &  \cdots \\
a & -a & 0 & 0 &  \cdots \\
0 & 0 & 0 & 0 & \cdots \\
\vdots & \vdots & \vdots & \vdots & \ddots
\end{array}
\right],
\end{equation}
where $a=A_{12}+A_{23}+A_{31}-A_{21}-A_{32}-A_{13}$.  Such combinations always vanish in the symmetric case, but are important in what follows.
Let
$$h(A)=\gcd(\{(A^{\pi \triangle})_{12}: \pi \in \mathcal{S}_n\}).$$
Of course, to avoid repetition it is enough to take one representative $\pi$ for each distinct pre-image of $(1,2,3)$.

\begin{lemma}
\label{lem:smallest-triangle}
Let $\Phi(A)$ be the primitivity lattice of $A$ and $h=h(A)$ as defined above.
Then $\mathrm{(1)}~h \mid b-c$ for every $(b,c) \in \Phi(A)$ and 
$\mathrm{(2)}~ (h,-h) \in \Phi(A)$.
%Then $h(A) =\frac{b-c}{2}$ or $b-c$.
\end{lemma}

\begin{proof}
Define matrices $E,E' \in \Z^{n \times n}$ by
$$
E_{ij} =
\begin{cases}
~1 & \text{if } (i,j)=(1,3) \text{ or } (2,4),\\
-1 & \text{if } (i,j)=(1,4) \text{ or } (2,3),\\
~0 & \text{otherwise}
\end{cases} 
\hspace{5mm}\text{and}
\hspace{5mm}
E'_{ij} = 
\begin{cases}
~1 & \text{if } (i,j)=(1,3), (3,2) \text{ or } (2,1),\\
-1 & \text{if } (i,j)=(1,2), (2,3) \text{ or } (3,1),\\
~0 & \text{otherwise.}
\end{cases} 
$$
For (1), it is easy to check that 
$(bE+cE^\top)+(bE+cE^\top)^{(123)}+(bE+cE^\top)^{(132)} = (b-c)E'$.
To prove (2), observe that 
\begin{equation}
\label{box-tri-relation}
A^\Box + A^{(123)\Box} + A^{(132)\Box} = A^\triangle+A^{\triangle(12)(34)}.
\end{equation}
Write $hE'=\sum_{\pi \in \Sym_n} t_\pi A^{\pi \triangle}$ for some coefficients $t_\pi \in \Z$.  Then
$$h(E-E^\top) = hE'+(hE')^{(12)(34)} = \sum_{\pi \in \mathcal{S}_n} t_\pi (A^{\pi \triangle} +A^{\pi \triangle (12)(34)}).$$
By \eqref{box-tri-relation}, this is a linear combination of copies of $A^{\pi \Box}$.  It follows that $(h,-h) \in \Phi(A)$.
\end{proof}

Some special properties hold for skew-symmetric matrices.
If $A+A^\top=O$, then each term in the gcd calculation of $h(A)$ has the form $2(A_{ij}+A_{jk}+A_{ki})$.  It follows that $h(A)$ is even.  In the special case that $h(A)=0$, we have 
$\Phi(A)=\{(0,0)\}$, which means the entries of $A$ satisfy $A_{ij} = f_i-f_j$ for some integers $f_1,\dots,f_n$.

Fix $A \in \Z^{n \times n}$ and put $h=h(A)$.
Define a triple of indices $\{i,j,k\} \subset [n]$ to be \emph{bad} if 
\begin{equation}
\label{eq:bad-test}
\nu_2(A_{ij}+A_{jk}+A_{ki}-A_{ji}-A_{kj}-A_{ik}) < 1+\nu_2(h).
\end{equation}
The idea behind the definition is as follows.  
Suppose we try to make local adjustments to some $B \in M_A$ to bring its off-diagonal entries closer to the average. Then, by parity, the six off-diagonal entries associated with a bad triple cannot be simultaneously corrected by a copy of $hE'$ on those indices, since $\nu_2(6h)=1+\nu_2(h)$.

For each $i \in [n]$, let $s_i=0$ or $1$ according as the number of bad triples containing $i$ is even or odd, respectively.  Put $\mathbf{s}(A)=(s_1,\dots,s_n) \in (\Z/2\Z)^n$. 
For a fixed $A \in \Z^{n \times n}$, a triple $\{i,j,k\}$ is bad in a sum of two matrices in $M_A$ if and only if it is bad in exactly one of the summands.  It follows that the map $A' \mapsto \mathbf{s}(A')$ from $M_A$ to $(\Z/2\Z)^n$ is a $\Z$-module homomorphism.  It also respects permutation; that is, $\mathbf{s}(A^\sigma) = \mathbf{s}^\sigma$.  
Therefore, a necessary condition on $\sum_\sigma c_\sigma A^\sigma = \left( \sum_\sigma c_\sigma \right) \overline{A}$ is that
\begin{equation}
\label{Aparity-nec}
\sum_{\sigma \in \Sym_n} c_\sigma \mathbf{s}^\sigma \equiv \mathbf{0} \pmod{2}.
\end{equation}
We also note that every pair of distinct indices $i<j$ are together in $n-2$ triples.  When $n$ is even, summing \eqref{eq:bad-test} over all index triples $i<j<k$ shows that $\sum_{i=1}^n s_i \equiv 0 \pmod{2}$ for any matrix $A$.

\subsection{The local lattice}

Let $A \in \Z^{n \times n}$ with off-diagonal row and column sums $d_i^+$ and $d_j^-$, respectively.  Let $G$ be a set of ($0$, $1$, or $2$) generators for the primitivity lattice $\Phi(A)$, as defined in Section~\ref{sec:prim-lattice}.  Let $(e_i^+,e_i^-)=(f_i^++f_i^--z^+,f_i^++f_i^--z^-)$.  
Let $s_i$ be the parity values for bad triples, as defined in Section~\ref{sec:parity}.

Define the \emph{local lattice} of $A$ as
\begin{equation}
\label{local-lattice}
\Lambda^*_A = \left\langle 
\begin{array}{c}
 \{(1,s_i,A_{ii},d_i^+,d_i^-,e_i^+,e_i^-):i=1,\dots,n\} \cup \{(0,2,0,0,0,0,0) \}  \\
\cup \{(0,0,0,0,0) \oplus \mathbf{u} : \mathbf{u} \in G\}
\end{array}
\right\rangle
\end{equation}
We have the following extension of Proposition~\ref{prop:G-local} lower bounding $\bal(A)$.

\begin{prop}
\label{prop:A-local}
Let $A \in \Z^{n \times n}$ with local lattice $\Lambda^*_A$.  Let $b$ denote the least positive multiple of 
$$\left( 1,0,\frac{\tr(A)}{n},\frac{\mathrm{sum}^*(A)}{n},\frac{\mathrm{sum}^*(A)}{n},\frac{\mathrm{sum}^*(A)}{n(n-1)},\frac{\mathrm{sum}^*(A)}{n(n-1)} \right)$$
belonging to $\Lambda^*_A$.  Then $b \mid \bal(A)$.
\end{prop}

\begin{proof}
Suppose $c_\sigma$ are coefficients such that $\sum_\sigma c_\sigma = \bal(A)$ and
\begin{equation}
\label{necessary-Acomb}
\sum_\sigma c_\sigma A^\sigma = \bal(A)\, \overline{A}.
\end{equation}
Let $C_i$ be the sum of $c_\sigma$ over all $\sigma \in \mathcal{S}_n$ such that $\sigma(1)=i$.
Then, examining the first row sum, column sum, and diagonal entry in \eqref{necessary-Acomb} gives
$$\sum_{i=1}^n C_i (1,A_{ii},d_i^+,d_i^-) = \bal(A) \left( 1, \frac{\tr(A)}{n}, \frac{\mathrm{sum}^*(A)}{n},\frac{\mathrm{sum}^*(A)}{n} \right).$$
From \eqref{Aparity-nec},
$$\sum_{i=1}^n C_i s_i = \sum_{\sigma \in \Sym_n} c_\sigma \mathbf{s}^\sigma_1 \equiv 0 \pmod{2}.$$
Put $\lam:=\bal(A) \, \frac{\mathrm{sum}^*(A)}{n(n-1)}$.
\comment{
The given combination implies something like
$$\sum_\sigma c_\sigma 
(f_i^+ + f_j^- - z^+,
f_j^+ + f_i^- - z^-) \equiv (\lam,\lam) \pmod{\Phi(A)}.$$
From \eqref{z-phi}, perhaps it doesn't matter if $\sigma$ reverses the order of $i,j$.
Classify according to perms that map $1 \mapsto i$ and $2 \mapsto j$.}
Let $C_{ij}$ be the sum of $c_\sigma$ over all $\sigma \in \mathcal{S}_n$ such that $\sigma(1)=i$ and $\sigma(2)=j$. 
%$$C_{ij}=\sum_{\sigma(1)=i,\sigma(2)=j} c_\sigma.$$
(This is zero if $i=j$.)
Next, consider the $(1,2)$ and $(2,1)$ entries in
\eqref{necessary-Acomb}.  On the right side, these off-diagonal entries are each $\lambda$.  From Lemma~\ref{lem:matrix-prim-dec} and \eqref{z-phi},
$$(A_{ij},A_{ji}) \equiv (F_{ij},F_{ji})+(Z_{ij},Z_{ji})\equiv
(f_i^+ + f_j^-,f_j^+ + f_i^-) -(z^+,z^-) \pmod{\Phi(A)}.$$ Therefore,
\begin{align*}
(\lam,\lam) &= \sum c_\sigma (A^\sigma_{12}, 
A^\sigma_{21})=
\sum_{i,j=1}^n C_{ij} (A_{ij},A_{ji})\\
&\equiv \sum_{i,j=1}^n C_{ij} (f_i^+ + f_j^- - z^+,
f_j^+ + f_i^- - z^-)\\
&= \sum_{i=1}^n C_{i} (f_i^+,f_i^-) + \sum_{j=1}^n C_j (f_j^-,f_j^+) - \sum_{i=1}^n C_i(z^+,z^-)\\
&= \sum_{i=1}^n C_{i} (e_i^+,e_i^-) \pmod{\Phi(A)}.\\
\end{align*}
We have shown that $\bal(A)$ times the given test vector belongs to $\Lambda^*(A)$, and so $b \mid \bal(A)$.
\end{proof}

\subsection{Sufficiency}

Here, we show the converse of Proposition~\ref{prop:A-local}, namely that the local
lattice $\Lambda^*_A$ is sufficient for computing the balancing index of $A$.
Similar to the case of graphs, we use the group ring elements $\alpha_r$.

Before the main result, we make some useful features about $F$ and $Z$ in the primitive decomposition. 

\begin{lemma}
\label{lem:FZ-features}
Suppose $A=B+F+Z$ in Lemma~\ref{lem:matrix-prim-dec}, and let $r \in [n]$.\\
{\rm (a)} For distinct indices $x,y>1$,
$F^{\alpha_r}_{xy} = f_r^+ + f_r^-$.\\
{\rm (a)} For distinct indices $x$ and $y$,
$(Z^{\alpha_r}_{xy},Z^{\alpha_r}_{yx})
\equiv -(z^+,z^-) \pmod{\Phi(A)}$
\end{lemma}

\begin{proof}
(a) Consider the case $r=1$, where 
$\alpha_1=\sum_{s=1}^n (1s) - (n-1) ()$. 
Since $x,y>1$, $F^{(1s)}_{xy}=f^+_x+f^-_y=F_{xy}$ unless $s=x$ or $s=y$.  These two terms of the sum not cancelling are $f^+_1$ and $f^-_1$, respectively.  Next we consider $r>1$, where
$\alpha_r=\sum_{s=1}^n (1sr) - (n-1) (1r)$.  If $x,y\neq r$, then just as above all terms cancel except for $f^+_r$ when $s=x$ and $f^-_r$ when $s=y$. But $(1sr)$ and $(1r)$ each map index $r$ to $1$, so the terms with $x=r$ or $y=r$ also cancel.  

(b) For any $\sigma \in \mathcal{S}_n$ and indices $i \neq j$, we have $(Z^\sigma_{ij},Z^\sigma_{ji})=-(z^+,z^-)$ or $-(z^-,z^+)$.  These are congruent $\pmod{\Phi(A)}$ by \eqref{z-phi}.  So, since the sum of coefficients in $\alpha_r$ equals $1$, we have  
$(Z^{\alpha_r}_{xy},Z^{\alpha_r}_{yx})
\equiv -(z^+,z^-) \pmod{\Phi(A)}$.
\end{proof}

\begin{thm}
\label{thm:matrices}
Let $A \in \Z^{n \times n}$ with local lattice $\Lambda^*_A$.  Let $C$ denote the least positive multiple of 
$$\left( 1,0,\frac{\tr(A)}{n},\frac{\mathrm{sum}^*(A)}{n},\frac{\mathrm{sum}^*(A)}{n},\frac{\mathrm{sum}^*(A)}{n(n-1)},\frac{\mathrm{sum}^*(A)}{n(n-1)} \right)$$
belonging to $\Lambda^*(A)$.  Then $\bal(A)=C$.
\end{thm}

\begin{proof}
Proposition~\ref{prop:A-local} shows $C \mid \bal(A)$, so it suffices to show the converse.
Put $\lam=C \frac{\mathrm{sum}^*(A)}{n(n-1)}.$
Let $C_i$ be integer coefficients summing to $C$ such that 
\begin{align}
\label{comb1}
\sum_{i=1}^n C_i (A_{ii},d_i^+,d_i^-) &= (C \, \tr(A)/n,\lam(n-1),\lam(n-1)),\\
\label{comb2}
\sum_{i=1}^n C_i s_i &\equiv 0 \pmod{2},\text{ and} \\
\label{comb3}
\sum_{i=1}^n C_i (e_i^+,e_i^-) &\equiv (\lam,\lam) \pmod{\Phi(A)}.
\end{align}
Put $K=\sum_{i=1}^n C_i A^{\alpha_i}$. From \eqref{comb1}, the diagonal entries of $K$ already match those of $C \overline{A}$.  Also from \eqref{comb1}, the row and column sums of $K$ agree with those of $C\overline{A}$.  
By \eqref{comb2}, any index belongs to an even number of bad triples.
We next show that $(K_{xy},K_{yx}) \equiv (\lam,\lam) \pmod{\Phi(A)}$.
Suppose $x,y>1$. Then from Lemma~\ref{lem:FZ-features} and \eqref{comb3},
\begin{align*}
(K_{xy},K_{yx}) \equiv \sum_{i=1}^n C_i (A^{\alpha_i}_{xy},A^{\alpha_i}_{yx}) &\equiv \sum_{i=1}^n C_i (F^{\alpha_i}_{xy},F^{\alpha_i}_{yx}) +
\sum_{i=1}^n C_i (Z^{\alpha_i}_{xy},Z^{\alpha_i}_{yx}) \\
&\equiv \sum_{i=1}^n C_i (e_i^+,e_i^-) 
\equiv (\lam,\lam) \pmod{\Phi(A)}.
\end{align*}

To handle the case when $x=1$ and $y>1$, we use the sum of the above over all $x \neq y$ to obtain
\begin{align*}
(K_{1y},K_{y1}) +(n-2)(\lam,\lam) &\equiv
\sum_{x \neq y} \sum_{i=1}^n C_i
(A_{xy}^{\alpha_i},A_{yx}^{\alpha_i})\\
&=\sum_{i=1}^n C_i \sum_{x \neq y} 
(A_{xy}^{\alpha_i},A_{yx}^{\alpha_i})\\
&=\sum_{i=1}^n C_i (d_i^+,d_i^-) =
\lam(n-1)(1,1).
\end{align*}
From this we have $(K_{1y},K_{y1})\equiv (\lam,\lam) \pmod{\Phi(A)}$ for each $y>1$,
and by a similar argument swapping rows with columns, 
$(K_{x1},K_{1x})\equiv (\lam,\lam) \pmod{\Phi(A)}$ for each $x>1$.

If $\Phi(A)=\{(0,0)\}$, we are already done.  If $\dim \Phi(A)=1$, we are either in the symmetric or skew-symmetric case, and take a basis vector $\mathbf{u}=(k,k)$ or $(k,-k)$.  In the generic case $\dim \Phi(A)=2$, we may choose a basis $\{\mathbf{u},\mathbf{v}\}$ where 
$\mathbf{u}=(u_1,u_2)$ and $\mathbf{v}=(u_2,u_1)$.

In what follows, we work with a $1$-norm relative to $\{\mathbf{u},\mathbf{v}\}$-coordinates, namely 
$\| \alpha \mathbf{u}+\beta \mathbf{v} \| := | \alpha | + | \beta |$, where the second term is omitted if $\dim \Phi(A)=1$.
Define $\mathrm{var}(K)=\sum_{i<j} \|(K_{ij},K_{ji})-(\lam,\lam)\|$.  
Suppose $K^*$ minimizes $\mathrm{var}(X)$ over all matrices $X \in M_A$ which have the same properties as $K$ described above.
For indices $x \neq y$, put $\delta_{xy}=K^*_{xy} - \lam$.  By \eqref{comb1}, we have $\sum_x \delta_{xy} = \sum_x \delta_{yx} = 0$ and by \eqref{comb3}, we have $(\delta_{xy},\delta_{yx}) \in \Phi(A)$.

Let us call a sequence of three pairs of indices $(x,w)$, $(x,y)$, $(z,y)$ $\mathbf{u}$-\emph{alternating} if the $\mathbf{u}$-coordinates of $(\delta_{xw},\delta_{wx})$, $(\delta_{xy},\delta_{yx})$, $(\delta_{zy},\delta_{yz})$ are nonzero and alternate in sign.  If $\mathrm{var}(K^*)>0$, there exist $x \neq y$ with $(\delta_{xy},\delta_{yx}) \neq (0,0)$.  
Suppose that the $\mathbf{u}$-coordinate is positive.  (The other cases are similar.)  From the above remarks, there exists some $z$ such that $(\delta_{zy},\delta_{yz})$ has negative $\mathbf{u}$-coordinate.  Likewise, there exists some $w$ such that $(\delta_{xw},\delta_{wx})$ has negative $\mathbf{u}$-coordinate.  Thus a $\mathbf{u}$ or $\mathbf{v}$-alternating sequence is guaranteed to exist, but we must consider cases according to whether it is `open' ($w \neq z$) or `closed' ($w=z$).

{\sc Case 1:}  $K^*$ has an open $\mathbf{u}$-alternating sequence.
Take an integer linear combination of copies of $A^\Box$, call it $U$, which vanishes except on its restriction to $\{x,y,z,w\}$, where it has the form
$$U[xzyw]=
\kbordermatrix{&z&x&y&w\\
x&\phantom{-}0 &\phantom{-}0 &\phantom{-}u_1 & -u_1  \\
z&\phantom{-}0 &\phantom{-}0 &-u_1 &\phantom{-}u_1   \\
y&\phantom{-}u_2 &-u_2 &\phantom{-}0 &\phantom{-}0   \\
w&-u_2 &\phantom{-}u_2 & \phantom{-}0 & \phantom{-}0 
}.$$
We obtain $\mathrm{var}(K^*-U)<\mathrm{var}(K^*)$ since at least three of the four modified terms have $\mathbf{u}$-coordinate reduced by one.  But $U \in M_A$ is formed with coefficient sum equal to $0$, has vanishing diagonal and line sums, and has off-diagonal pairs in $\Phi(A)$.  This is a contradiction to the choice of $K^*$.

{\sc Case 2:}  Every $\mathbf{u}$- or $\mathbf{v}$-alternating sequence in $K^*$ is closed.
By assumption, each of the ordered pairs
$(\delta_{xy},\delta_{yx})$, $-(\delta_{zy},\delta_{yz})$, $-(\delta_{xz},\delta_{zx})$ has  $\mathbf{u}$-coordinate of the same sign, say positive.  Suppose for the moment that the $\mathbf{v}$-coordinates were also positive.  Then since 
$\mathbf{u}+\mathbf{v}$ is a constant vector, the analysis in the proof of Theorem~\ref{thm:multigraphs} (the symmetric case) would produce an open alternating sequence, contradicting our assumption for this case. 
Suppose, via Lemma~\ref{box-tri-relation}(1), that $\mathbf{u}-\mathbf{v}$ is a positive multiple of $(h,-h)$, where $h=h(A)$ is as defined in Section~\ref{sec:parity}.  (If negative, we simply let $-h$ take the role of $h$ below.)

Take an integer linear combination of copies of $A^{\pi \triangle}$, call it $U'$, 
which vanishes except on its restriction to $\{x,y,z\}$,
$$U'[xyz]=
\kbordermatrix{&x&y&z\\
x&\phantom{-}0 & \phantom{-}h & -h  \\
y& -h & \phantom{-}0 &\phantom{-}h  \\
z&\phantom{-}h & -h & \phantom{-}0 }.
$$
Consider $K^*-U'$.  From our choices above, we have 
$\mathrm{var}(K^*-U') \le \mathrm{var}(K^*)$.  We claim the inequality is strict.  
Recall that the number of bad triples in $K^*$ containing each of $x,y$ is even.  If $\{x,y,z\}$ were bad, we would again contradict that there are no open alternating sequences.  Therefore,
$$2h \mid \delta_{xy}+\delta_{yz}+\delta_{zx}-\delta_{xz}-\delta_{zy}-\delta_{yx}.$$
It follows that $\|(\delta_{xy},\delta_{yx})-(h,-h)\|+\|(\delta_{xz},\delta_{zx})+(h,-h)\|+
\|(\delta_{zy},\delta_{yz})+(h,-h)\|$ is indeed less than $\|(\delta_{xy},\delta_{yx})\|+\|(\delta_{xz},\delta_{zx})\|+\|(\delta_{zy},\delta_{yz})\|.$
As before, this contradicts choice of $K^*$.  We must have $\mathrm{var}(K^*) = 0$, and thus $K^*=C \overline{A}$ as required.
\end{proof}

As in the case of multigraphs, the proof suggests an algorithm that iteratively reduces the deviation from $C \overline{A}$ to zero.  We give two examples.

\begin{ex}
Consider $$A=
\begin{bmatrix}
0 & 1 & 4 & 1 \\
0 & 0 & 0 & 1 \\
0 & 3 & 0 & 3 \\
1 & 5 & 1 & 0
\end{bmatrix}.$$
A routine calculation gives 
$\Phi(A)=\langle (4,1),(1,4) \rangle$ and $h(A)=3$.  The vectors of row/column sums 
are $(6,1,6,7)$ and $(1,9,5,5)$.  The lists of $e_i^+$ and $e_j^-$ can each be taken as $(1,0,3,1)$.
The bad triples of indices are $\{1,2,4\}$ and $\{2,3,4\}$, giving parity vector $\mathbf{s}=(1,0,1,0)$.  Even without a lattice computation, we see that $\overline{A}=\frac{5}{3} (J-I)$, so it is necessary that $3 \mid \bal(A)$.  Using the lattice $\Lambda^*_A$, Theorem~\ref{thm:matrices} predicts $\bal(A)=3$, with coefficients $C_1=2$, $C_2=2$, $C_3=-8$, $C_4=7$.  From these, the first step is to compute $K=2A^{\alpha_1}+2A^{\alpha_2}-
8A^{\alpha_3}+7A^{\alpha_4}$,
which has balanced line sums, 
satisfies the parity condition on bad triples, and has off-diagonal pairs $(K_{ij},K_{ji}) \equiv 3(\overline{A}_{ij},\overline{A}_{ji}) = (5,5) \equiv (0,0) \pmod{\Phi(A)}$.

We describe a sequence of adjustments to 
$$[\delta_{xy}]=K-3\overline{A}=
\left[
\begin{array}{rrrr}
0 & \color[rgb]{0,0,1} \mathbf{-24} & -8 & \color[rgb]{0,0,1} \mathbf{32}\\
-6 & 0 & 11 & -5\\
\color[rgb]{1,0,0} \mathbf{43} & -16 & 0 & \color[rgb]{0,0,1} \mathbf{-27}\\
\color[rgb]{1,0,0} \mathbf{-37}  & \color[rgb]{1,0,0} \mathbf{40} & -3 & 0
\end{array}
\right]$$
using copies of
$U=A^\Box=
\mbox{\scriptsize
$\left[
\begin{array}{rrrr}
0 & 0 & 4 & -4 \\
0 & 0 & -4 & 4 \\
1 & -1 & 0 & 0 \\
-1 & 1 & 0 & 0
\end{array}
\right]$}$
and 
$U'=
A^{(34)\triangle}=
\mbox{\scriptsize
$\left[
\begin{array}{rrrr}
0 & -3 & 3 & 0 \\
3 & 0 & -3 & 0 \\
-3 & 3 & 0 & 0 \\
0 & 0 & 0 & 0 
\end{array}
\right]$}$.\\
The highlighted entries above have large and alternating deviation from $\lam$. If we subtract $6U^{(13)(24)}$ and add $6U^{(23)}$, this aligns a correction of $\pm 18$ or $\pm 24$ on these entries. The result of these adjustments is
$$
\left[
\begin{array}{rrrr}
0 & 0 & -14 & 14\\
0 & 0 & 11 & -11\\
19 & -16 & 0 & -3\\
-19 & 16 & 3 & 0
\end{array}\right].$$
It is straightforward to check that $(-14,19) \in \Phi(A)$, so a further adjustment leaves
$$
\left[
\begin{array}{rrrr}
0 & 0 & 0 & 0\\
0 & 0 & -3 & 3\\
0 & 3 & 0 & -3\\
0 & -3 & 3 & 0
\end{array}\right].$$
Subtracting a copy of $U'$ reduces this to the zero matrix, finishing the verification that $\bal(A)=3$.
\end{ex}

\begin{ex}
Consider $$A=
\left[
\begin{array}{rrrrr}
0 & 1 & 2 & 2 & 1\\
2 & 0 & 0 & 0 & 1\\
2 & 1 & 0 & 0 & 0\\
0 & 1 & 2 & 0 & 0\\
2 & 1 & 1 & 1 & 0
\end{array}\right].$$
We have $\Phi(A)=\langle (1,1),(1,-1) \rangle$ and $h(A)=2$.  After a calculation,
Theorem~\ref{thm:matrices} says $\bal(A)=1$ and $\Lambda_A^*$ gives the locally balanced matrix $K$ below.  
After adjusting the marked entries by copies of $\pm 2E$ or $\pm 4E$, one obtains the  matrix shown at right, which happens to be symmetric.
$$K=
\left[
\begin{array}{rrrrr}
0 & 0 & \color[rgb]{.4,0,1} \mathbf{-2} & \color[rgb]{.4,0,1} \mathbf{6} & 0\\
0 & 0 & \color[rgb]{.4,0,1} \mathbf{5} & \color[rgb]{.4,0,1} \mathbf{-3} & 2\\
2 & \color[rgb]{0,.6,0} \mathbf{-1} & 0 & 0 & \color[rgb]{0,.6,0} \mathbf{3}\\
\color[rgb]{1,.6,0} \mathbf{6} & \color[rgb]{0,.6,0} \mathbf{3} & \color[rgb]{1,.6,0} \mathbf{-4} & 0 & \color[rgb]{0,.6,0} \mathbf{-1}\\
\color[rgb]{1,.6,0} \mathbf{-4} & 2 & \color[rgb]{1,.6,0} \mathbf{5} & 1 & 0
\end{array}\right] \hspace{5mm} \longrightarrow
\hspace{5mm}
\left[
\begin{array}{rrrrr}
0 & 0 & 2 & 2 & 0\\
0 & 0 & 1 & 1 & 2\\
2 & 1 & 0 & 0 & 1\\
2 & 1 & 0 & 0 & 1\\
0 & 2 & 1 & 1 & 0
\end{array}\right].$$
The deviation from $\overline{A}=J-I$ is represented by a signed `bowtie' graph, as shown in Figure~\ref{fig:bowtie}.
Being a sum of two copies of $C_4^*$, this deviation can be reduced to the zero matrix.
Note that, although $K$ has bad triples, none remain in the bowtie since the underlying matrix is symmetric.
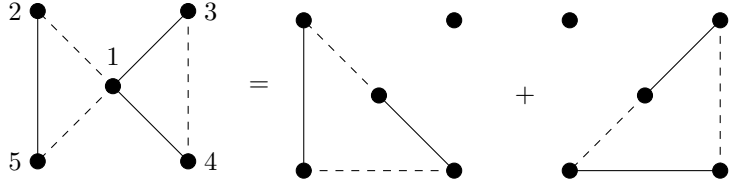
\begin{figure}[htbp]
\begin{center}
\begin{tikzpicture}
\draw[dashed] (0,1)--(-1,0);
\draw[dashed] (0,1)--(-1,2);
\draw (-1,0)--(-1,2);
\draw (0,1)--(1,0);
\draw (0,1)--(1,2);
\draw[dashed] (1,0)--(1,2);
\filldraw (-1,0) circle [radius=0.1];
\filldraw (-1,2) circle [radius=0.1];
\filldraw (0,1) circle [radius=0.1];
\filldraw (1,0) circle [radius=0.1];
\filldraw (1,2) circle [radius=0.1];
\node at (-1.3,0) {$5$};
\node at (-1.3,2) {$2$};
\node at (0,1.4) {$1$};
\node at (1.3,0) {$4$};
\node at (1.3,2) {$3$};
\node at (2,1) {$=~~$};
\end{tikzpicture}
~\begin{tikzpicture}
\draw[dashed] (0,1)--(-1,2);
\draw (-1,0)--(-1,2);
\draw (0,1)--(1,0);
\draw[dashed] (-1,0)--(1,0);
\filldraw (-1,0) circle [radius=0.1];
\filldraw (-1,2) circle [radius=0.1];
\filldraw (0,1) circle [radius=0.1];
\filldraw (1,0) circle [radius=0.1];
\filldraw (1,2) circle [radius=0.1];
\node at (2,1) {$+~~$};
\end{tikzpicture}
~\begin{tikzpicture}
\draw[dashed] (0,1)--(-1,0);
\draw (0,1)--(1,2);
\draw[dashed] (1,0)--(1,2);
\draw (-1,0)--(1,0);
\filldraw (-1,0) circle [radius=0.1];
\filldraw (-1,2) circle [radius=0.1];
\filldraw (0,1) circle [radius=0.1];
\filldraw (1,0) circle [radius=0.1];
\filldraw (1,2) circle [radius=0.1];
\end{tikzpicture}
\caption{A compound symmetric adjustment; negative edges dashed}
\label{fig:bowtie}
\end{center}
\end{figure}
\end{ex}

A minor restatement of Theorem~\ref{thm:matrices} classifies the completely symmetric matrices in $M_A$.

\begin{cor}
We have $\mu I + \lam J \in M_A$ if and only if 
$(C,0,\lam+\mu,\lam(n-1),\lam(n-1),\lam,\lam) \in \Lambda^*_A$.
\end{cor}

We also remark that our lattice condition yields a general upper bound on $\bal(A)$.

\begin{cor}
\label{cor:bound}
Let $A \in \Z^{n \times n}$.  Then $\bal(A) \mid n(n-1)$.
\end{cor}

\begin{proof}
Consider $\mathbf{a}=(n-1) \sum_{i=1}^n (1,s_i,A_{ii},d_i^+,d_i^-,e_i^+,e_i^-)$.  The first component is $n(n-1)$.  The second is even, using the remark at the end of Section~\ref{sec:parity}.  The third equals $(n-1) \tr(A)$.  The fourth and fifth equal $(n-1) \mathrm{sum}^*(A)$.  It remains to consider the last two coordinates.
Using Lemma~\ref{lem:matrix-prim-dec} and \eqref{z-phi}, 
$$(n-1) \sum_{i=1}^n (e_i^+,e_i^-) =
(n-1) \sum_{i=1}^n (f^+_i+f^-_i) (1,1) - n(n-1)(z^+,z^-) \equiv \mathrm{sum}^*(A) (1,1) \pmod{\Phi(A)}.$$
It follows that $\mathbf{a}$ equals $n(n-1)$ times the test vector in Theorem~\ref{thm:matrices}, and thus $\bal(A) \mid n(n-1)$.
\end{proof}

For prime powers $n$, there is a na\"ive proof of Corollary~\ref{cor:bound}.  Simply index entries with a finite field and sum over the group of affine transformations.  It is not clear if such a direct argument for this bound exists for general $n$.

\subsection{Noteworthy cases}

In each of the symmetric and skew-symmetric cases, we can simplify our lattice condition.   
A symmetric matrix $A$ with index of primitivity $k$ has $\Phi(A)=\langle (k,k) \rangle$.
We also have $d_i^+=d_i^-=d_i$ and can replace $e_i^+,e_i^-$ with potentials $f_i$.  This gives the multigraph lattice with an additional coordinate to track diagonal entries.

\begin{prop}
\label{prop:sym}
Suppose $A \in \Z^{n \times n}$ is symmetric with off-diagonal row sums $d_i$, potentials $f_i$, and index of primitivity $k$.
Then $\bal(A)$ is the least positive multiple of $(1,\frac{\tr(A)}{n},\frac{\mathrm{sum}^*(A)}{n},\frac{\mathrm{sum}^*(A)}{n(n-1)})$ belonging to 
$$\langle \{(1,A_{ii},d_i,2f_i):i=1,\dots,n\} \cup \{(0,0,0,k)\} \rangle \subseteq \Z^4.$$
\end{prop}

In the skew-symmetric case, we have $A_{ii}=0$, $d_i^+ = -d_i^-$, and can take $e_i^+ = -e_i^-$, $f_i=0$.  Also, $\Phi(A)=\langle (k,-k) \rangle$ for some integer $k \ge 0$, which is analogous to the index of primitivity in the symmetric case.  These specializations result in the following simplified lattice condition.
\begin{prop}
\label{prop:skew}
Suppose $A \in \Z^{n \times n}$ is skew-symmetric with off-diagonal row sums $d_i^+$ and parity values $s_i$.  Then $\bal(A)$ is the least positive multiple of $(1,0,0,0)$ belonging to 
$$\langle \{(1,s_i,d^+_i,e^+_i):i=1,\dots,n\} \cup
\{(0,2,0,0),(0,0,0,k)\} \rangle \subseteq \Z^4.$$
\end{prop}

A special case of skew-symmetric matrices has a graph-theoretic interpretation.
A \emph{tournament} $T$ is a digraph in which, for any two vertices $x \neq y$, exactly one of $(x,y)$ and $(y,x)$ is an arc of $T$.   The sequence of vertex outdegrees is also known as the \emph{score sequence} of $T$.  The adjacency matrix of $T$ is a matrix $A \in \{0,1\}^{n \times n}$ satisfying $A+A^\top = J-I$. That is, $2A-J+I$ is skew-symmetric with entries $\pm 1$.

Every tournament on at least four vertices has a transitive triangle, i.e. three distinct arcs of the form $x \rightarrow y \rightarrow z$, $x \rightarrow z$ as well as four arcs of the form $w \rightarrow x \rightarrow y \rightarrow z$, $w \rightarrow z$.  Let $A$ again denote the adjacency matrix of $T$.  We can use these configurations to obtain $\Z$-linear combinations of $A^{\pi \triangle}$ realizing $h(A)=1$ and $\Z$-linear combinations of $A^{\pi \Box}$ realizing $\Phi(A) = \langle (1,-1) \rangle$.
That is, we have $k=1$ in Proposition~\ref{prop:skew}.

\begin{prop}
\label{prop:T}
Let $T$ be a tournament on $n \ge 2$ vertices with score sequence $d_1^+,\dots,d_n^+$.  Then 
$\bal(T)=2$ if $T$ is regular, and otherwise
$$\bal(T) = \frac{2q}{\gcd(q,n-1-2r)},$$
where $q=\gcd(d_i^+-d_j^+: 1 \le i,j \le n)$ and $d_i^+ \equiv r \pmod{q}$ for each $i=1,\dots,n$.
\end{prop}

\begin{proof}
Note first that $\bal(T)$ is an even integer, since the average off-diagonal entry in its adjacency matrix equals $\frac{1}{2}$.  Consider the special case when 
$T$ is regular, so that $n$ is odd and each $d_i^+= \frac{n-1}{2}$.  The lattice condition becomes vacuous, and it follows that $\bal(T)=2$.

Now suppose $T$ is not regular, and let $q=\gcd(d_i^+-d_j^+: 1 \le i,j \le n)$.  
Consider the least positive integer $b$ such that $b(1,n-1) \in \langle (1,2d_i^+) : i = 1,\dots,n \rangle$.  Note that $b \mid n$ since $2\sum d_i^+ = n(n-1)$.  By considering entries and row sums it is easy to see that $b \mid \bal(T)$.  Conversely, from the above discussion, we need not consider other coordinates in the local lattice, and hence $\bal(T)=b$.

We now compute a closed formula for $b$.  By Lemma~\ref{lem:2d-basis}, our lattice $\langle (1,2d_i^+)\rangle$ is generated by 
$(1,2r)$, $(0,2q)$.  Then, by Lemma~\ref{lem:2d-closedform}, we have $b=2q/\gcd(q,n-1-2r)$.  
\end{proof}

Another noteworthy class of digraphs are the $1$-regular ones, which also correspond to permutation matrices. 

\begin{prop}
\label{prop:P}
Suppose $\pi \in \mathcal{S}_n$ has exactly $t$ fixed points.  Let $P_\pi$ be the permutation matrix corresponding to $\pi$.  Then
\begin{equation}
\label{eq:perm-matrix}
\bal(P_\pi) = \frac{n(n-1)}{\gcd(n-t,(n-1)t)},
\end{equation}
unless $\pi$ has cycle type $(4)$ or $(3,3)$, in which case $\bal(P_\pi)$ is twice the value in \eqref{eq:perm-matrix}.
\end{prop}

\begin{proof}[Proof sketch.]
The result for $n \le 6$ can be verified by direct computation for each of the integer partitions of $n$.  Suppose $n \ge 7$, and let $D$ be the digraph associated with $\pi$.  If $D$ is symmetric (that is, if $\pi$ has no cycle of length at least $3$), then  $\bal(P_\pi)$ can be computed via the underlying simple graph consisting of $(n-t)/2$ disjoint edges on $n$ vertices.  In this case, our formula \eqref{eq:perm-matrix} agrees with \eqref{balG-primitive}.

Suppose now that $\pi$
contains a cycle of length at least $3$. It is straightforward to check that there exist four vertices of $D$ which induce exactly one arc.  For instance, one can choose all four vertices on a cycle of length at least $7$, three of the vertices on a cycle of length $5$ or $6$, or an arc of some shorter cycle along with two other non-adjacent vertices. This implies the primitivity lattice $\Phi(P_\pi)$ is simply $\Z^2$, and the minimum triangle adjusment is $h(P_\pi)=1$.  After examining a few cases, the vector $\mathbf{s}$ is either $(0,\dots,0)$ or $(1,\dots,1)$ according to the parity of $n-t$. However, if $n-t$ is odd, then \eqref{eq:perm-matrix} returns an even value.  Thus the local lattice for $P_\pi$ can be reduced to just three dimensions and only three distinct generators: $(1,1,0)$ for fixed points, $(1,0,1)$ for non-fixed points, and $(0,0,1)$ for the primitivity lattice.  The test vector is $(1,\frac{t}{n},\frac{n-t}{n(n-1)})$, and it suffices to take the least positive multiple which belongs to $\Z^3$.  This corresponds with the right side of \eqref{eq:perm-matrix}, as required.
\end{proof}

To close, we remark that the $3 \times 3$ case, examined earlier in Section~\ref{sec:3x3}, can be reformulated with a similar appearance as the result for $n \ge 4$.  Here, $\Phi(A)$ disappears, as there is no room to construct a combination such as $A^\Box$.  
As a caution to the reader, we have kept diagonal elements separate in this section, but they were included in the line sums of Section~\ref{sec:3x3} and the summary result below.

\begin{prop}
Suppose $A \in \Z^{3 \times 3}$ is symmetric with line sums $d_i^+$, $d_j^-$.
Then $\bal(A)$ is the least positive multiple of $(1,0,\frac{1}{3}\tr(A),\frac{1}{3}\mathrm{sum}(A),\frac{1}{3}\mathrm{sum}(A))$ belonging to 
$$\langle \{(1,\epsilon,A_{ii},d_i^+,d_i^-):i=1,2,3\} \cup \{(0,2,0,0,0)\} \rangle,$$
where $\epsilon=0$ if $A^\triangle=O$, and otherwise $\epsilon=1$.
\end{prop}

\begin{proof}
Let $b$ be the least multiple in the statement.  First, we claim $b \mid \bal(A)$. Working from \eqref{necessary-Acomb} and following the start of the proof of Proposition~\ref{prop:A-local}, we see $\bal(A) (1,\frac{1}{3}\tr(A),\frac{1}{3}\mathrm{sum}(A),\frac{1}{3}\mathrm{sum}(A))$ is a $\Z$-linear combination of $(1,A_{ii},d_i^+,d_i^-)$.
The second coordinate forces $b$ to be even in the event that  $A_{12}+A_{23}+A_{31} \neq A_{21}+A_{32}+A_{13}$.  But Lemma~\ref{lem:parity3x3} also implies $\bal(A)$ is even in this case. 

Next, we claim $\bal(A) \mid b$.  Again, by Lemma~\ref{lem:parity3x3}, we have that $b$ is even if $\bal(A)$ is even.  Suppose $3 \mid \bal(A)$.  Then, by Proposition~\ref{prop:3x3}, either $\mathrm{rank}(\Theta(A)) = 3$, or at least one of $A_{ii}$, $d_i^+$, $d_i^-$, $i=1,2,3$, is not a ternary triple.  Any of these conditions implies that there is no integer solution to $\sum_{i=1}^3 x_i (1,A_{ii},d_i^+,d_i^-) = b(1,\frac{1}{3}\tr(A),\frac{1}{3}\mathrm{sum}(A),\frac{1}{3}\mathrm{sum}(A))$ unless $3 \mid b$.
\end{proof}

\section{Conclusion}

We have studied the module generated by copies of a given matrix (graph), classifying the
completely symmetric matrices (respectively, multiples of complete graphs) belonging to this module.
To this end, we defined a balancing index for matrices and graphs.  A general method for computing balancing index is described, recovering a variety of existing results on signed graph designs and extending to new settings.  The unifying feature in all cases is that a suitably defined local lattice (much smaller than the `full' module) governs which completely symmetric matrices or multiples of complete graphs are present in the relevant permutation module.

In Table~\ref{table:progression}, we summarize various results on this problem, ordered by complexity of the local lattice.

\begin{table}[htbp]
\begin{tabular}{lcl}
\hline
setting & dim $\Lambda^*$ & reference \\
\hline
graphs with two isolated vertices & 1 & \cite[Proposition 4]{WW}; also \cite{Wilson75} \\
graphs with loops and two isolates  & 2 & \cite[Lemma 2.3]{DM} \\
digraphs with two isolates  & 2 & \cite[Proposition 4]{Wilson75} \\
primitive multigraphs & 2 & \cite[Theorem 13]{WW} \\
general multigraphs & 3 & Theorem~\ref{thm:multigraphs}\\
symmetric matrices & 4 & Proposition~\ref{prop:sym} \\
skew-symmetric matrices & 4 & Proposition~\ref{prop:skew} \\
general matrices & 7 & Theorem~\ref{thm:matrices}\\
\hline
\end{tabular}
\medskip
\caption{Progression of signed balancing results with local lattice dimensions}
\label{table:progression}
\end{table}

One topic for future work is a full description of the relevant module, say by Smith forms.  Since much of the ground work already exists in \cite{WW}, the next step is to introduce potential and diagonal-entry data into these calculations.  It is also natural to consider the homology involved in higher-rank extensions, say to hypergraphs or tensors.  Our problem extended to hypergraphs is connected with the concept of a `null design', \cite{GJ,MS,WilsonSHD}.

In balancing a multigraph or matrix, we have allowed all permutations to be used.  A natural next question is to consider a subgroup $\Gamma \le \Sym_n$ of allowed permutations and aim to balance $A$ with least positive coefficient sum $\bal_\Gamma(A)$.  In contrast with Corollary~\ref{cor:bound}, we found that there exist $5 \times 5$ matrices whose balancing index over the alternating group is $60$.  When $\Gamma$ is $2$-transitive, a completely symmetric matrix always results from the sum of all $A^\sigma$, $\sigma \in \Gamma$.  Even if $\Gamma$ is not $2$-transitive, one might classify the matrices $A$ which can be balanced under $\Gamma$.

A related notion is balancing in a block-diagonal sense.  To offer a concrete instance of this, consider an $r$-partite graph $G$ on vertex set $[n] \times [r]$ with vertex partition $\{[n] \times \{1\}, \dots, [n] \times \{r\}\}$.   If we allow those permutations $\Sym_n \wr \Sym_r$ which stabilize the partition, one can ask for the least positive sum of coefficients producing a $\lambda$-fold complete $r$-partite graph.
%; in terms of adjacency matrices, the target is $
%\lam J_n \otimes (J_r-I_r)$. 
The special case when $G$ consists of $n$ disjoint cliques $K_r$ transverse to the partition gives a signed analog of resolvable transversal designs.

Another variation arises from a restriction that coefficients $c_\sigma$ be nonnegative integers.  This can be viewed as an extension of the problem of computing Frobenius numbers, and also includes the existence question for balanced incomplete block designs.  For these reasons, the problem is expected to be hard in general.  Still, there are likely accessible and interesting questions in this direction.  Applications of designs, say to numerical integration or the design of statistical experiments, may invite extensions to a multigraph or matrix setting and serve as additional motivation.

%\section*{Declarations}
%Declarations of interest: none.

\section*{Acknowledgments}

Research of Peter Dukes is supported by NSERC grant RGPIN 2017--03891.

\end{document}